\catcode`\@=11

\font\tenmsa=msam10
\font\sevenmsa=msam7
\font\fivemsa=msam5
\font\tenmsb=msbm10
\font\sevenmsb=msbm7
\font\fivemsb=msbm5
\newfam\msafam
\newfam\msbfam
\textfont\msafam=\tenmsa  \scriptfont\msafam=\sevenmsa
  \scriptscriptfont\msafam=\fivemsa
\textfont\msbfam=\tenmsb  \scriptfont\msbfam=\sevenmsb
  \scriptscriptfont\msbfam=\fivemsb

\def\hexnumber@#1{\ifnum#1<10 \number#1\else
 \ifnum#1=10 A\else\ifnum#1=11 B\else\ifnum#1=12 C\else
 \ifnum#1=13 D\else\ifnum#1=14 E\else\ifnum#1=15 F\fi\fi\fi\fi\fi\fi\fi}

\def\msa@{\hexnumber@\msafam}
\def\msb@{\hexnumber@\msbfam}
\mathchardef\boxdot="2\msa@00
\mathchardef\boxplus="2\msa@01
\mathchardef\boxtimes="2\msa@02
\mathchardef\square="0\msa@03
\mathchardef\blacksquare="0\msa@04
\mathchardef\centerdot="2\msa@05
\mathchardef\lozenge="0\msa@06
\mathchardef\blacklozenge="0\msa@07
\mathchardef\circlearrowright="3\msa@08
\mathchardef\circlearrowleft="3\msa@09
\mathchardef\rightleftharpoons="3\msa@0A
\mathchardef\leftrightharpoons="3\msa@0B
\mathchardef\boxminus="2\msa@0C
\mathchardef\Vdash="3\msa@0D
\mathchardef\Vvdash="3\msa@0E
\mathchardef\vDash="3\msa@0F
\mathchardef\twoheadrightarrow="3\msa@10
\mathchardef\twoheadleftarrow="3\msa@11
\mathchardef\leftleftarrows="3\msa@12
\mathchardef\rightrightarrows="3\msa@13
\mathchardef\upuparrows="3\msa@14
\mathchardef\downdownarrows="3\msa@15
\mathchardef\upharpoonright="3\msa@16

\mathchardef\downharpoonright="3\msa@17
\mathchardef\upharpoonleft="3\msa@18
\mathchardef\downharpoonleft="3\msa@19
\mathchardef\rightarrowtail="3\msa@1A
\mathchardef\leftarrowtail="3\msa@1B
\mathchardef\leftrightarrows="3\msa@1C
\mathchardef\rightleftarrows="3\msa@1D
\mathchardef\Lsh="3\msa@1E
\mathchardef\Rsh="3\msa@1F
\mathchardef\rightsquigarrow="3\msa@20
\mathchardef\leftrightsquigarrow="3\msa@21
\mathchardef\looparrowleft="3\msa@22
\mathchardef\looparrowright="3\msa@23
\mathchardef\circeq="3\msa@24
\mathchardef\succsim="3\msa@25
\mathchardef\gtrsim="3\msa@26
\mathchardef\gtrapprox="3\msa@27
\mathchardef\multimap="3\msa@28
\mathchardef\therefore="3\msa@29
\mathchardef\because="3\msa@2A
\mathchardef\doteqdot="3\msa@2B

\mathchardef\triangleq="3\msa@2C
\mathchardef\precsim="3\msa@2D
\mathchardef\lesssim="3\msa@2E
\mathchardef\lessapprox="3\msa@2F
\mathchardef\eqslantless="3\msa@30
\mathchardef\eqslantgtr="3\msa@31
\mathchardef\curlyeqprec="3\msa@32
\mathchardef\curlyeqsucc="3\msa@33
\mathchardef\preccurlyeq="3\msa@34
\mathchardef\leqq="3\msa@35
\mathchardef\leqslant="3\msa@36
\mathchardef\lessgtr="3\msa@37
\mathchardef\backprime="0\msa@38
\mathchardef\risingdotseq="3\msa@3A
\mathchardef\fallingdotseq="3\msa@3B
\mathchardef\succcurlyeq="3\msa@3C
\mathchardef\geqq="3\msa@3D
\mathchardef\geqslant="3\msa@3E
\mathchardef\gtrless="3\msa@3F
\mathchardef\sqsubset="3\msa@40
\mathchardef\sqsupset="3\msa@41
\mathchardef\trianglerighteq="3\msa@44
\mathchardef\trianglelefteq="3\msa@45
\mathchardef\bigstar="0\msa@46
\mathchardef\between="3\msa@47
\mathchardef\blacktriangledown="0\msa@48
\mathchardef\blacktriangleright="3\msa@49
\mathchardef\blacktriangleleft="3\msa@4A
\mathchardef\blacktriangle="0\msa@4E
\mathchardef\triangledown="0\msa@4F
\mathchardef\eqcirc="3\msa@50
\mathchardef\lesseqgtr="3\msa@51
\mathchardef\gtreqless="3\msa@52
\mathchardef\lesseqqgtr="3\msa@53
\mathchardef\gtreqqless="3\msa@54
\mathchardef\Rrightarrow="3\msa@56
\mathchardef\Lleftarrow="3\msa@57
\mathchardef\veebar="2\msa@59
\mathchardef\barwedge="2\msa@5A
\mathchardef\doublebarwedge="2\msa@5B
\mathchardef\angle="0\msa@5C
\mathchardef\measuredangle="0\msa@5D
\mathchardef\sphericalangle="0\msa@5E
\mathchardef\varpropto="3\msa@5F
\mathchardef\smallsmile="3\msa@60
\mathchardef\smallfrown="3\msa@61
\mathchardef\Subset="3\msa@62
\mathchardef\Supset="3\msa@63
\mathchardef\Cup="2\msa@64

\mathchardef\Cap="2\msa@65

\mathchardef\curlywedge="2\msa@66
\mathchardef\curlyvee="2\msa@67
\mathchardef\leftthreetimes="2\msa@68
\mathchardef\rightthreetimes="2\msa@69
\mathchardef\subseteqq="3\msa@6A
\mathchardef\supseteqq="3\msa@6B
\mathchardef\bumpeq="3\msa@6C
\mathchardef\Bumpeq="3\msa@6D
\mathchardef\lll="3\msa@6E

\mathchardef\ggg="3\msa@6F

\mathchardef\circledS="0\msa@73
\mathchardef\pitchfork="3\msa@74
\mathchardef\dotplus="2\msa@75
\mathchardef\backsim="3\msa@76
\mathchardef\backsimeq="3\msa@77
\mathchardef\complement="0\msa@7B
\mathchardef\intercal="2\msa@7C
\mathchardef\circledcirc="2\msa@7D
\mathchardef\circledast="2\msa@7E
\mathchardef\circleddash="2\msa@7F
\def\ulcorner{\delimiter"4\msa@70\msa@70 }
\def\urcorner{\delimiter"5\msa@71\msa@71 }
\def\llcorner{\delimiter"4\msa@78\msa@78 }
\def\lrcorner{\delimiter"5\msa@79\msa@79 }
\def\yen{\mathhexbox\msa@55 }
\def\checkmark{\mathhexbox\msa@58 }
\def\circledR{\mathhexbox\msa@72 }
\def\maltese{\mathhexbox\msa@7A }
\mathchardef\lvertneqq="3\msb@00
\mathchardef\gvertneqq="3\msb@01
\mathchardef\nleq="3\msb@02
\mathchardef\ngeq="3\msb@03
\mathchardef\nless="3\msb@04
\mathchardef\ngtr="3\msb@05
\mathchardef\nprec="3\msb@06
\mathchardef\nsucc="3\msb@07
\mathchardef\lneqq="3\msb@08
\mathchardef\gneqq="3\msb@09
\mathchardef\nleqslant="3\msb@0A
\mathchardef\ngeqslant="3\msb@0B
\mathchardef\lneq="3\msb@0C
\mathchardef\gneq="3\msb@0D
\mathchardef\npreceq="3\msb@0E
\mathchardef\nsucceq="3\msb@0F
\mathchardef\precnsim="3\msb@10
\mathchardef\succnsim="3\msb@11
\mathchardef\lnsim="3\msb@12
\mathchardef\gnsim="3\msb@13
\mathchardef\nleqq="3\msb@14
\mathchardef\ngeqq="3\msb@15
\mathchardef\precneqq="3\msb@16
\mathchardef\succneqq="3\msb@17
\mathchardef\precnapprox="3\msb@18
\mathchardef\succnapprox="3\msb@19
\mathchardef\lnapprox="3\msb@1A
\mathchardef\gnapprox="3\msb@1B
\mathchardef\nsim="3\msb@1C
\mathchardef\napprox="3\msb@1D
\mathchardef\nsubseteqq="3\msb@22
\mathchardef\nsupseteqq="3\msb@23
\mathchardef\subsetneqq="3\msb@24
\mathchardef\supsetneqq="3\msb@25
\mathchardef\subsetneq="3\msb@28
\mathchardef\supsetneq="3\msb@29
\mathchardef\nsubseteq="3\msb@2A
\mathchardef\nsupseteq="3\msb@2B
\mathchardef\nparallel="3\msb@2C
\mathchardef\nmid="3\msb@2D
\mathchardef\nshortmid="3\msb@2E
\mathchardef\nshortparallel="3\msb@2F
\mathchardef\nvdash="3\msb@30
\mathchardef\nVdash="3\msb@31
\mathchardef\nvDash="3\msb@32
\mathchardef\nVDash="3\msb@33
\mathchardef\ntrianglerighteq="3\msb@34
\mathchardef\ntrianglelefteq="3\msb@35
\mathchardef\ntriangleleft="3\msb@36
\mathchardef\ntriangleright="3\msb@37
\mathchardef\nleftarrow="3\msb@38
\mathchardef\nrightarrow="3\msb@39
\mathchardef\nLeftarrow="3\msb@3A
\mathchardef\nRightarrow="3\msb@3B
\mathchardef\nLeftrightarrow="3\msb@3C
\mathchardef\nleftrightarrow="3\msb@3D
\mathchardef\divideontimes="2\msb@3E
\mathchardef\varnothing="0\msb@3F
\mathchardef\nexists="0\msb@40
\mathchardef\mho="0\msb@66
\mathchardef\thorn="0\msb@67
\mathchardef\beth="0\msb@69
\mathchardef\gimel="0\msb@6A
\mathchardef\daleth="0\msb@6B
\mathchardef\lessdot="3\msb@6C
\mathchardef\gtrdot="3\msb@6D
\mathchardef\ltimes="2\msb@6E
\mathchardef\rtimes="2\msb@6F
\mathchardef\shortmid="3\msb@70
\mathchardef\shortparallel="3\msb@71
\mathchardef\smallsetminus="2\msb@72
\mathchardef\thicksim="3\msb@73
\mathchardef\thickapprox="3\msb@74
\mathchardef\approxeq="3\msb@75
\mathchardef\succapprox="3\msb@76
\mathchardef\precapprox="3\msb@77
\mathchardef\curvearrowleft="3\msb@78
\mathchardef\curvearrowright="3\msb@79
\mathchardef\digamma="0\msb@7A
\mathchardef\varkappa="0\msb@7B
\mathchardef\hslash="0\msb@7D
\mathchardef\hbar="0\msb@7E
\mathchardef\backepsilon="3\msb@7F
\def\Bbb{\ifmmode\let\next\Bbb@\else
 \def\next{\errmessage{Use \string\Bbb\space only in math mode}}\fi\next}
\def\Bbb@#1{{\Bbb@@{#1}}}
\def\Bbb@@#1{\fam\msbfam#1}

\catcode`\@=\active

\def\sw#1{{\sb{(#1)}}}
\def\sco#1{{\sp{(\bar #1)}}} 
\def\su#1{{\sp{(#1)}}} 
\def\d{{\rm d}}
\def\proof{{\sl Proof.}\ }
\def\endproof{\hbox{$\sqcup$}\llap{\hbox{$\sqcap$}}}
\def\tens{\mathop{\otimes}}

\def\CL{{\Lambda}}
\def\CM{{\cal M}}
\def\CN{{\cal N}} 
\def\o{{}_{(1)}}
\def\t{{}_{(2)}}
\def\Bo{{}_{\und{(1)}}}
\def\Bt{{}_{\und{(2)}}}
\def\th{{}_{(3)}}
\def\Bth{{}_{\und{(3)}}}
\def\<{{\langle}}
\def\>{{\rangle}}
\def\und#1{{\underline{#1}}}
\def\ra{{\triangleleft}}
\def\la{{\triangleright}} 
\def\id{{\rm id}} 
\def\eps{\epsilon}
\def\span{{\rm span}}
\def\q2{{q^{-2}}}
\def\bicross{{\blacktriangleright\!\!\!\triangleleft}}
\def\note#1{{}}
\def\eqn#1#2{\begin{equation}#2\label{#1}\end{equation}}
\def\qbinom#1#2#3{\left(\begin{array}{c}#1\\#2\end{array}\right)\sb#3}
\def\Z{{\Bbb Z}}

\documentstyle[12pt,epsf]{article}

\newtheorem{prop}{Proposition}[section]

\newtheorem{lemma}[prop]{Lemma}

\newtheorem{df}[prop]{Definition}
\newtheorem{ex}[prop]{Example}

\begin{document}
\baselineskip 22pt

{\ }\qquad\qquad \hskip 4.3in DAMTP/95-74

\vspace{.2in}

\begin{center} {\LARGE COALGEBRA GAUGE THEORY\footnote{Research
supported by the EPSRC grant GR/K02244}}
\\ \baselineskip 13pt{\ }
{\ }\\
Tomasz Brzezi\'nski and Shahn Majid\footnote{Royal Society University
Research Fellow and
Fellow of
 Pembroke College, Cambridge. On leave 1995 + 1996 at the Department of
 Mathematics, Harvard University, Cambridge MA02138, USA}\\
{\ }\\
 Department of Applied Mathematics \& Theoretical Physics\\
University of Cambridge, Cambridge CB3 9EW\\
\end{center}
\begin{center}
December 1995  -- revised January 1996
\end{center}

\vspace{10pt}
\begin{quote}\baselineskip 13pt
\noindent{\bf Abstract} We develop a generalised gauge theory in which
the role
of gauge group is played by a coalgebra and the role of principal
bundle by an
algebra. The theory provides a unifying point of view which includes quantum
group gauge theory, embeddable quantum homogeneous spaces and braided group
gauge theory, the latter being introduced now by these means. Examples
include
ones in which the gauge groups are the braided line and the quantum  plane.

\bigskip
\noindent Keywords:  gauge theory -- coalgebra -- quantum group -- braided
group -- quantum plane -- connection -- bicrossproduct -- factorisation

\end{quote}
\baselineskip 19pt

\section{Introduction}

In a recent paper \cite{Brz:hom} it was shown by the first author that
a generalisation of the quantum group principal bundles introduced in
\cite{BrzMa:gau} is needed if one wants to include certain
`embeddable' quantum   homogeneous spaces, such as the full family of
quantum two-spheres of Podle\'s   \cite{Pod:sph}. A one-parameter
specialisation of this family was used in   \cite{BrzMa:gau}  in
construction of the $q$-monopole, but the general members of the
family do not have the required canonical fibering. The required
generalised notion of quantum principal bundles proposed in \cite{Brz:hom},
also termed a $C$-Galois
extension, consists of an algebra $P$, a coalgebra $C$ with a
distinguished element $e$ and a right action of $P$ on $P\otimes C$
satisfying certain conditions. In the present paper we develop a version of
such  `coalgebra principal bundles' based  on a map $\psi:C\tens P\to P\tens C$
and $e\in C$, and giving now a theory of connections on  them.

Another motivation for the paper is the search for a generalisation of gauge
theory powerful enough to include braided
groups\cite{Ma:bra}\cite{Ma:bg}\cite{Ma:introp} as the gauge group. Although
not quantum groups, braided groups do have at least a coalgebra and hence can
be covered in our theory.  We describe the main elements of such
a braided group gauge theory as  arising in this way. This is a first step
towards a theory of braided-Lie algebra valued gauge fields, Chern-Simmons and
Yang-Mills actions, to be considered elsewhere.

As well as providing a unifying point of view which includes our previous
quantum group gauge theory\cite{BrzMa:gau}, the theory of embeddable
homogeneous spaces\cite{Brz:hom} and braided group gauge theory, our
coalgebra
gauge theory has its own characteristic properties. In particular, the axioms
obeyed by $\psi$ involve the algebra and coalgebra in a symmetrical way,
opening up the possibility of an interesting self-duality of the
construction.
This becomes manifest when we are given a character $\kappa$ on $P$; then we
have also the possibility of a dual `algebra gauge theory', corresponding in
the finite-dimensional case to a coalgebra gauge theory with the coalgebra
$P^*$, principal bundle total space $C^*$ and the structure map $\psi^*$.
This is a new phenomenon which is not possible within the realm of ordinary
(non-Abelian) gauge theory. Moreover, the axioms obeyed by $\psi$
correspond in
the finite-dimensional case to the factorisation of an algebra into $P^{\rm
op}C^*$, which is a common situation\cite{Ma:phy}. Indeed, all bicrossproduct
quantum groups\cite{Ma:phy} provide a dual pair of examples.

Finally, we note that some steps towards a theory of fibrations based on
algebra factorisations has appeared independently in \cite{CapMic:non},
including topological considerations which may be useful in further work.
However, we really need the present coalgebra treatment for our
infinite-dimensional algebraic examples, for our treatment of differential
calculus and in order to include  quantum and braided group gauge theories.
We demonstrate the various stages of our formalism on some concrete examples
based on the braided line and quantum plane.
\vspace{12pt}

{\sc Preliminaries.}  All vector spaces are taken over a field
$k$ of generic characteristic.
$C$ denotes a
coalgebra with the coproduct $\Delta: C\to C\otimes C$ and the
counit $\epsilon: C\to k$ which satisfy the standard axioms.
For the coproduct we use the Sweedler notation
$$
\Delta c = c\sw 1\otimes c\sw 2,\qquad \Delta\sp 2 c = (\Delta
\otimes \id)\circ\Delta c = c\sw 1\otimes c\sw 2
\otimes c\sw 3,\quad {\rm etc.},
$$
where $c\in C$, and the summation sign and the
indices are suppressed.

A vector space
$P$ is a right $C$-comodule if there exists
a map $\Delta\sb R : P\to P\otimes C$, such that
$(\Delta\sb R\otimes \id)\circ\Delta\sb R =
(\id\otimes\Delta)\circ\Delta\sb R$, and
$(\id\otimes \epsilon)\circ\Delta\sb R = \id$.
For $\Delta\sb R$ we use the explicit notation
$$
\Delta\sb R u = u\sco 0\otimes u\sco 1 ,
$$
where $u\in P$ and $u\sco 0\otimes c\sco 1\in P\tens C$ (summation
understood).
For $e\in C$, we denote by $P\sb e\sp{coC}$
the vector subspace of $P$ of all elements $u\in P$
such that $\Delta\sb R u = u\otimes e$.

$H$ denotes a Hopf algebra with product $\cdot : H\otimes H
\to H$, unit $1$, coproduct $\Delta: H\to H\otimes H$,
counit $\epsilon: H\to k$
and antipode $S: H\to H$. We use Sweedler's sigma notation
as before. Similarly as for a coalgebra we can
define right  $H$-comodules. We say that
a right  $H$-comodule $P$ is a
right  $H$-comodule algebra if $P$
is an algebra and $\Delta\sb R$
is an algebra map.

If $P$ is an algebra then by $\Omega\sp n P$ we denote a $P$-bimodule
of universal $n$-forms on $P$, i.e.
$$
\Omega\sp n P = \{ \omega\in P\sp{\otimes n+1}:\;\; \forall i\in
\{1,\ldots, n\},\;\; \cdot\sb i \omega =0\},
$$
where $\cdot\sb i$ denotes a multiplication in $P$ acting on the $i$ and
$i+1$ factors in $P\sp{\otimes n+1}$. $\Omega\sp n(P)$ denotes a bimodule of
$n$-forms on $P$ obtained from
$\Omega\sp n P$ as an appropriate quotient.

Finally, we recall the definition of a $C$-Galois extension from
\cite{Brz:hom}.
This consists  of $C$ a coalgebra, $P$ an algebra, $\Delta_R:P\to P\tens C$ a
right coaction of $C$ on $P$ and $\rho :P\otimes C\otimes
P\to P\otimes C$ a right action of $P$ on $P\otimes C$, such that
\eqn{C-Galois.df}{\Delta\sb R\circ\cdot = \rho\circ (\Delta\sb R\otimes
\id).}
We also require an element  $e\in C$ such that  $\rho(u\otimes e, v) =uv\sco
0\otimes v\sco 1$, for any $u,v\in P$. Then $M=P\sb e\sp{co C}$ is an algebra
and we say that
$P(M,C,\rho,e)$ is a {\em $C$-Galois extension} or a {\em quantum
$\rho$-principal bundle over $M$} if the canonical map
$\chi\sb M:P\otimes\sb MP\to P\otimes C$, $\chi\sb M :u\otimes\sb
Mv\mapsto uv\sco 0\otimes v\sco 1$ is a bijection.

In \cite{Brz:hom} we provided two main examples of $C$-Galois extensions.
Firstly, it is shown that if $H$ is a Hopf algebra and  $P$ a right
$H$-comodule algebra (so $M=P\sp{coH}$ is an algebra) and if the map $\chi\sb
M: P\otimes\sb
MP\to P\otimes H$
as defined above is a bijection then the Hopf-Galois extension $M\subset P$
(i.e. the {\em quantum
principal bundle} $P(M,H)$ with universal differential structure) is a
$C$-Galois extension with $C=H$, $e=1$ and $\rho(u\otimes c ,v) = uv\sco
0\otimes cv\sco 1$.

Secondly, let $H$ be a Hopf algebra, $C$ a coalgebra and
$\pi: H\to C$ a coalgebra projection. Then $H$ is a right $C$-comodule
with a coaction $\Delta\sb R =  (\id\otimes\pi)\circ\Delta$.  Denote
$e =\pi(1)\in C$ and define $M=H\sb e\sp{coC}$ as before. Assume
that $\ker\pi$ is a minimal right ideal in $H$ such that
$\{ u - \epsilon(u) ; u\in M\} \subset \ker\pi$. Then \cite{Brz:hom}
\eqn{C-Galois.hom}{
\rho(u\otimes c, v) = uv\sw 1\otimes \pi (wv\sw 2),}
for any $u,v\in H$, $c\in C$ and $w\in\pi\sp{-1}(c)$ makes $H(M,C,\rho, e)$ a
$C$-Galois extension or
 quantum $\rho$-principal bundle associated to an {\em embeddable quantum
homogeneous space}. Some concrete examples are also given in \cite{Brz:hom}
(see also \cite{DijKoo:hom}).

In addition to these gauge-theoretic preliminaries, we will also discuss
examples based on the theory of braided
groups\cite{Ma:bra}\cite{Ma:bg}\cite{Ma:introp} and the theory of
bicrossproduct and double cross product and Hopf algebras
\cite{Ma:phy}\cite{Ma:rem}, due to the second
author. Chapters~6.2,7.2,9 and~10
of the text \cite{Ma:book} contain full details on these topics.

\section{Coalgebra $\psi$-Principal Bundles}
In this paper we will be dealing with a particular class of $C$-Galois
extensions or generalised quantum principle bundles. This class is more
tractable than the most general case in \cite{Brz:hom}, allowing us to
develop
a gauge theory for it in the next section.  Yet, it is general enough to
include all our main examples of interest. Our data is the following:

\begin{df} We say that  a coalgebra $C$  and an algebra $P$ are {\em
entwined}
if there is a map  $\psi :C\otimes P\to
P\otimes C$ such that
\eqn{ent.A}{\psi\circ(\id\tens
\cdot)=(\cdot\tens\id)\circ\psi_{23}\circ\psi_{12},\quad \psi(c\tens
1)=1\tens   c,\quad \forall c\in C}
\eqn{ent.B}{(\id\tens\Delta)\circ\psi=\psi_{12}\circ\psi_{23}
\circ(\Delta\tens\id),\quad (\id\tens\eps)\circ\psi=\eps\tens\id,}
where $\cdot$ denotes multiplication in $P$, and
$\psi_{23}=\id\tens\psi$ etc.   Explicitly, we require that the
following diagrams commute:
\begin{equation}
\begin{picture}(455,125)(20,1)
\put(40,100){\makebox(0,0){$C\otimes P\otimes P$}}
\put(70,100){\vector(1,0){100}}
\put(107,105){$\id\otimes \cdot$}
\put(190,100){\makebox(0,0){$C\otimes P$}}
\put(190,90){\vector(0,-1){50}}
\put(200,65){$\psi$}
\put(190,35){\makebox(0,0){$P\otimes C$}}
\put(40,90){\vector(0,-1){50}}
\put(50,65){$\psi\otimes \id$}
\put(40,35){\makebox(0,0){$P\otimes C\otimes P$}}
\put(40,27){\vector(2,-1){46}}
\put(85,16){\makebox(0,0){$\id\otimes\psi$}}
\put(115,1){\makebox(0,0){$P\otimes P\otimes C$}}
\put(142,3){\vector(2,1){46}}
\put(188,12){\makebox(0,0){$\cdot\otimes \id$}}
\put(250,100){\makebox(0,0){$C\otimes k$}}
\put(280,100){\vector(1,0){100}}
\put(317,105){$\id\otimes \eta$}
\put(400,100){\makebox(0,0){$C\otimes P$}}
\put(400,90){\vector(0,-1){50}}
\put(410,65){$\psi$}
\put(400,35){\makebox(0,0){$P\otimes C$}}
\put(248,90){\line(0,-1){50}}
\put(252,90){\line(0,-1){50}}
\put(250,35){\makebox(0,0){$k\otimes C$}}
\put(278,35){\vector(1,0){102}}
\put(320,40){\makebox(0,0){$\eta\otimes \id$}}
\end{picture}
\label{diag.A}
\end{equation}
\begin{equation}
\begin{picture}(455,135)(20,1)
\put(40,100){\makebox(0,0){$P\otimes C\otimes C$}}
\put(170,100){\vector(-1,0){100}}
\put(107,105){$\id\otimes \Delta$}
\put(190,100){\makebox(0,0){$P\otimes C$}}
\put(190,40){\vector(0,1){50}}
\put(200,65){$\psi$}
\put(190,35){\makebox(0,0){$C\otimes P$}}
\put(40,40){\vector(0,1){50}}
\put(50,65){$\psi\otimes \id$}
\put(40,35){\makebox(0,0){$C\otimes P\otimes C$}}
\put(86,4){\vector(-2,1){46}}
\put(85,16){\makebox(0,0){$\id\otimes\psi$}}
\put(115,1){\makebox(0,0){$C\otimes C\otimes P$}}
\put(188,28){\vector(-2,-1){46}}
\put(188,12){\makebox(0,0){$\Delta\otimes \id$}}
\put(250,100){\makebox(0,0){$P\otimes k$}}
\put(380,100){\vector(-1,0){100}}
\put(317,105){$\id\otimes \epsilon$}
\put(400,100){\makebox(0,0){$P\otimes C$}}
\put(400,40){\vector(0,1){50}}
\put(410,65){$\psi$}
\put(400,35){\makebox(0,0){$C\otimes P$}}
\put(248,90){\line(0,-1){50}}
\put(252,90){\line(0,-1){50}}
\put(250,35){\makebox(0,0){$k\otimes P$}}
\put(380,35){\vector(-1,0){102}}
\put(320,40){\makebox(0,0){$\epsilon\otimes \id$}}
\end{picture}
\label{diag.B}
\end{equation}
\label{entwined.df}
\end{df}

In the finite-dimensional case this is exactly equivalent by partial
dualisation to the requirement that $\tilde\psi:C^*\tens P^{\rm op}\to P^{\rm
op}\tens C^*$ is an algebra factorisation structure (which is part of the
theory of Hopf algebra double cross products\cite{Ma:phy}). This is made
precise at the end of the section, where it provides a natural way to obtain
examples of such $\psi$.

\begin{prop}
Let $C,P$ be entwined by $\psi$. For every group-like element $e\in C$
we have  the following:

1. For any positive $n$, $P\sp{\otimes n}$ is a right $C$-comodule
with the coaction
$\Delta\sb R\sp n = \psi\sb{nn+1}\circ\psi\sb{n-1n}\circ\ldots\circ\psi\sb
{12}\circ (\eta\sb C\otimes \id\sp n) \equiv \overleftarrow{\psi}\sp
n\circ(\eta\sb C\otimes \id\sp n)$,  where $\eta\sb C: k\to C$,
$\alpha\mapsto \alpha e$.

2. The coaction $\Delta\sb R\sp n$ restricts to the coaction on
$\Omega\sp n P$.

3. $M=P\sb e\sp{coC} = \{u\in P\; ;\;\; \Delta\sb R\sp 1 u = u\otimes
e\}$ is a subalgebra of $P$.

4. The linear map $\chi\sb M : P\otimes\sb M P \to P\otimes C$,
$u\otimes\sb M v\mapsto u\psi(e\otimes v)$ is well-defined. If
$\chi\sb M$
is a bijection we say that we have a $\psi$-principal bundle $P(M,C,\psi,e)$.

5. $P(M,C,\psi,e)$ is an example of a $C$-Galois extension with $\rho =
(\cdot\otimes \id)\circ\psi\sb {23}$.
\label{bundle.prop}
\end{prop}
\proof
We write $\psi(c\otimes u) = \sum\sb \alpha u\sb \alpha\otimes c\sp
\alpha$ and  henceforth we omit the summation sign. In this notation,
the conditions
(\ref{diag.A}) and (\ref{diag.B})  are
\begin{equation}
(uv)\sb \alpha\otimes c\sp \alpha = u\sb \alpha v\sb \beta\otimes
c\sp{\alpha\beta},
\qquad 1\sb \alpha\otimes c\sp \alpha = 1\otimes c,
\label{ind.A}
\end{equation}
\begin{equation}
u\sb \alpha\otimes {c\sp \alpha}\o \otimes {c\sp \alpha}\t =
u\sb{\alpha\beta}\otimes {c\o}\sp \beta\otimes {c\t}\sp \alpha, \qquad
\epsilon(c\sp \alpha)u\sb \alpha = \epsilon(c) u,
\label{ind.B}
\end{equation}
for all $u,v\in P$ and $c\in C$.

1. The map $\Delta\sp n\sb R$ is given explicitly by
$$
\Delta\sb R\sp n (u\sp 1\otimes \ldots\otimes u\sp n) = u\sp 1\sb {\alpha\sb
1}\otimes\ldots\otimes u\sp n\sb{\alpha\sb n}
\otimes e\sp{\alpha\sb 1\ldots \alpha\sb n}.
$$
Hence
\begin{eqnarray*}
(\Delta\sp n\sb R\otimes \id)\Delta\sb R\sp n
(u\sp 1\otimes \ldots\otimes u\sp n) &=&  u\sp 1\sb {\alpha\sb 1 \beta\sb 1}
\otimes\ldots\otimes u\sp n\sb{\alpha\sb n\beta\sb n}\otimes e\sp{\beta\sb
1\ldots \beta\sb n}
\otimes e\sp{\alpha\sb 1\ldots \alpha\sb n}\\
&=&   u\sp 1\sb {\alpha\sb 1 \beta\sb 1}
\otimes\ldots\otimes u\sp n\sb{\alpha\sb n\beta\sb n}\otimes {e\sw
1}\sp{\beta\sb
1\ldots \beta\sb n}
\otimes {e\sw 2}\sp{\alpha\sb 1\ldots \alpha\sb n}\\
&=&   u\sp 1\sb {\alpha\sb 1 }\otimes u\sp
2\sb{\alpha \sb 2 \beta \sb 2}
\otimes\ldots\otimes u\sp n\sb{\alpha\sb n\beta\sb n}\otimes
 {e\sp{\alpha\sb 1}\sw 1}{}\sp{\beta\sb 2\ldots \beta\sb n}\otimes
 {e\sp{\alpha\sb 1}\sw 2}{}\sp{\alpha\sb 2\ldots \alpha\sb n}\\
&=& ... =   u\sp 1\sb {\alpha\sb 1}
\otimes\ldots\otimes u\sp n\sb{\alpha\sb n}\otimes e\sp{\alpha\sb 1\ldots
\alpha\sb n}\sw 1
\otimes e\sp{\alpha\sb 1\ldots \alpha\sb n}\sw 2\\
&=& (\id\sp n\otimes \Delta)\Delta\sb R\sp n
(u\sp 1\otimes \ldots\otimes u\sp n),
\end{eqnarray*}
where we used the group-like property of $e$ to derive the second
equality and then we used the condition (\ref{diag.B}) $n$ times to
obtain the penultimate one. We also have
\begin{eqnarray*}
(\id\sp n\otimes\epsilon)\Delta\sb R\sp n
(u\sp 1\otimes \ldots\otimes u\sp n) &=& u\sp 1\sb {\alpha\sb
1}\otimes\ldots\otimes u\sp n\sb{\alpha\sb n}
\epsilon( e\sp{\alpha\sb 1\ldots \alpha\sb n}) \\
&=& u\sp 1\sb {\alpha\sb 1}\otimes\ldots\otimes u\sp n
\epsilon( e\sp{\alpha\sb 1\ldots \alpha\sb{n-1}}) \\
&=&\ldots = \epsilon(e)u\sp 1\otimes \ldots\otimes u\sp n\\
&=& u\sp 1\otimes \ldots\otimes u\sp n,
\end{eqnarray*}
where we have first used the condition (\ref{diag.B}) $n$-times and then the
group-like property of $e$. Hence $\Delta\sb R\sp n$ is a coaction.

2. If $\sum\sb i u\sp i\otimes v\sp i \in\ker \cdot$
then $(\cdot\otimes \id)\sum\sb i \Delta\sb R\sp 2 (u\sp i\otimes v\sp i) =
\sum\sb{i} u\sb \alpha\sp i v\sp i\sb \beta \otimes e\sp{\alpha\beta}  =
\sum\sb{i} (u\sp i v\sp i)\sb \alpha\otimes e\sp \alpha =0$, using
(\ref{ind.A}). Hence the coaction preserves $\Omega\sp 1 P$. Similarly for
$\Omega^nP$.

3. Here $M=\{u\in P| u_\alpha \tens e^\alpha=u\tens e\}$ and if $u,v\in M$ then
$(uv)_\alpha\tens e^\alpha=u_\alpha v_\beta\tens
e^{\alpha\beta}=uv_\beta\tens e^\beta=uv\tens e$ as well, using
(\ref{ind.A}).

4. It is easy to see that $\chi_M$ is well-defined as a map from
$P\tens_M P$.
Thus, if $x\in M$ we have $\chi_M(u,xv)=u(xv)_\alpha\tens e^\alpha=u x_\alpha
v_\beta\tens e^{\alpha\beta}=ux v_\beta\tens e^\beta=\chi_M(ux,v)$, using
(\ref{ind.A}).

5. Parts 3 and 4  also follow from the theory of $C$-Galois
extensions. Indeed,
we can define $\rho$ as stated and verify that
\begin{eqnarray*}
\rho\circ\Delta\sb
R\sp 1(u\otimes v)&=&(\cdot\otimes \id)\circ\psi\sb
{23}\circ(\Delta\sb R\sp 1\otimes \id)(u\otimes
v)
 =  (\cdot\otimes \id)\circ\psi\sb {23}(u\sb \alpha\otimes e\sp \alpha
\otimes v) \\
& = & u\sb \alpha v\sb \beta\otimes e\sp{\alpha\beta}
\stackrel{(\ref{ind.A})}{=}
   (uv)\sb \alpha\otimes e\sp \alpha
= \Delta\sb R\sp 1(uv).
\end{eqnarray*}
Furthermore,
$ \rho(u\otimes e, v) = u\psi(e,v) = u\Delta\sb R\sp 1v$ by definition of the
coaction $\Delta\sb R\sp 1$. Also from (\ref{ind.A}), it is easy to see that
$\rho$ is a right action, as required for a $C$-Galois extension.
\endproof
\begin{ex}
Let $H$ be a Hopf algebra and $P$ be a right $H$-comodule
algebra. The linear map $\psi :H\otimes P\to P\otimes H$ defined by
$\psi: c\otimes u\to u\sco 0\otimes cu\sco 1$ entwines $H,P$. Therefore a
quantum group principal
bundle $P(M,H)$ with universal differential structure as in \cite{BrzMa:gau}
 is a $\psi$-principal bundle
$P(M,H,\psi,1)$.
\label{quantum.ex}
\end{ex}
\proof
For any $c\in H$ and $u\in P$ we have $u\sb \alpha\otimes c\sp \alpha
=u\sco 0\otimes cu\sco  1$. Clearly $1\sb \alpha\otimes c\sp \alpha =
1\otimes c$. We compute
\[
u\sb \alpha v\sb \beta\otimes c\sp{\alpha\beta} =  u\sco 0v\sb
\beta\otimes (cu\sco 1)\sp \beta =u\sco 0v\sco 0\otimes cu\sco 1v\sco 1=
(uv)\sco 0
\otimes c(uv)\sco 1 = (uv)\sb \alpha\otimes c\sp \alpha,
\]
hence the condition (\ref{diag.A}) is satisfied. Furthermore,
$\epsilon(c\sb \alpha)u\sp \alpha = \epsilon(cu\sco 1)u\sco 0 =
\epsilon(c)u$ and
\begin{eqnarray*}
u\sb{\alpha\beta}\otimes c\sw 1\sp \beta\otimes c\sw 2\sp \alpha & =
& u\sco 0\sb \beta\otimes c\sw 1\sp \beta\otimes c\sw 2u\sco 1 = u\sco 0\sco
0\otimes
c\sw 1u\sco 0\sco 1\otimes c\sw 2u\sco 1 \\
& = & u\sco 0\otimes (cu\sco 1)\sw 1\otimes (cu\sco 1)\sw 2 =
u\sb \alpha\otimes c\sp \alpha\sw 1\otimes c\sb \alpha\sw 2,
\end{eqnarray*}
so that the condition (\ref{diag.B}) is also satisfied. Clearly the
induced coaction in Proposition~\ref{bundle.prop} coincides with the given
coaction of $H$.
\endproof

We can easily replace $H$ here by one of the braided groups introduced in
\cite{Ma:bra}\cite{Ma:bg}. To be concrete, we suppose that our braided
group $B$ lives in a $k$-linear
braided category with well-behaved direct sums, such as that of modules over
a quasitriangular Hopf algebra or comodules over a dual-quasitriangular Hopf
algebra.
This background quantum group does not enter directly into the braided group
formulae
but rather via the braiding $\Psi$ which it induces between any
objects in the
category.
We refer to \cite{Ma:introp} for an introduction to the theory and for
further
details.
In particular, a right braided $B$-module algebra $P$ means a coaction $P\to
P\und\tens B$
in the category which is an algebra homomorphism to the braided tensor
product
algebra\cite{Ma:bra}
\eqn{btens}{(u\tens b)(v\tens c)=u\Psi(b\tens v)c.}
The coproduct $\und\Delta:B\to B\und \tens B$ of a braided group is itself a
homomorphism to
such a braided tensor product.

\begin{ex}
Let $B$ be a braided group with braiding $\Psi$ and $P$ a right braided
$B$-comodule
algebra. The linear map $\psi :B\otimes P\to P\otimes B$ defined by
$\psi: c\otimes u\to \Psi(c\tens u\sco 0)u\sco 1$ entwines $B,P$. If the
induced map $\chi_M$ is
a bijection we say that the associated $\psi$-principal bundle
$P(M,B,\psi,1)$ is a {\em braided group principal bundle}, and denote it by
$P(M,B,\Psi)$.
\label{braided.ex}
\end{ex}
\proof This is best done diagrammatically by the technique introduced in
\cite{Ma:bra}. Thus, we write $\Psi=\epsfbox{braid.eps}$ and products by
$\cdot=\epsfbox{prodfrag.eps}$. We denote coactions and coproducts by
$\epsfbox{deltafrag.eps}$. The proof of the main part of
(\ref{diag.A}) is then   the diagram:
\[ \epsfbox{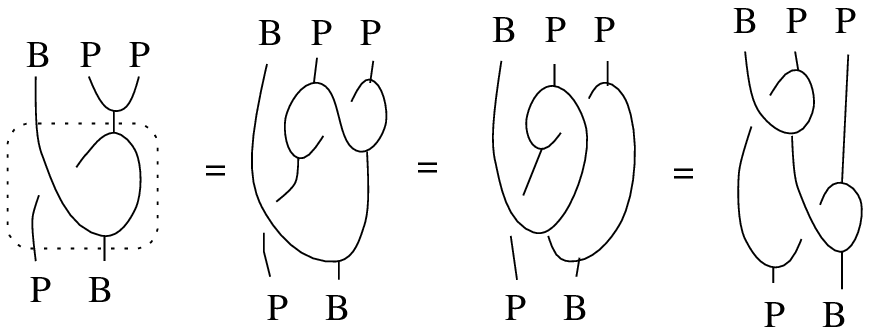}\]
where box is $\psi$ as stated, in diagrammatic form. The first
equality is the  assumed homomorphism  property of the braided coaction
$\epsfbox{deltafrag.eps}$. The second equality is associativity of the
product   in $B$, and the third is functoriality of the braiding,
which we use to push  the diagram into the right form. The minor
condition is immediate from the
axioms of a braided comodule algebra and the  properties of the unit map
$\eta:\underline{1}\to P$. Here $\underline{ 1}$ denotes the trivial
object for
our tensor product and necessarily commutes with the braiding in an
obvious way
(such that $\underline{1}$ is denoted consistently by omission). For
the proof
of (\ref{diag.B}) we ask the reader to reflect the diagram in a mirror
about a
horizontal axis (i.e. view it up-side-down and from behind) and then reverse
all braid crossings (restoring them all to $\epsfbox{braid.eps}$). The result
is the diagrammatic proof for the main part of (\ref{diag.B}) if we
relabel the
product of $P$ as the coproduct of $B$ and relabel the   product of
$B$  as the
right coaction of $B$ on $P$. The minor part of (\ref{diag.B}) is immediate
from properties of the braided counit.  \endproof

\begin{ex}
Let $H$ be a Hopf algebra and   $\pi :H\to C$ a coalgebra surjection. Then
$\psi :C\otimes H\to H\otimes C$ defined by $\psi(c\otimes u) =u\sw 1\otimes
\pi(vu\sw 2)$
entwines $C,H$, where $u\in H$, $c\in C$ and $v\in\pi\sp{-1}(c)$. If
$\ker\pi$
is a minimal right ideal containing $\{u-\eps(u)|u\in M\}$ then we have a
$\psi$-principal bundle
$H(M,C,\psi,\pi(1))$ in the setting of Proposition~\ref{bundle.prop}, denoted
$H(M,C,\psi,\pi)$.
Hence the generalised bundles over embeddable quantum homogeneous spaces
in \cite{Brz:hom} are examples of
$\psi$-principal
bundles.
\label{embeddable.ex}
\end{ex}
\proof
In this case $u\sb \alpha\otimes c\sp \alpha = u\sw 1\otimes\pi(wu\sw
2)$, for any $u\in H$, $c\in C$ and $w\in\pi\sp{-1}(c)$.  Clearly
$1\sb \alpha\otimes c\sp \alpha = 1\otimes c$. We compute
$$
u\sb \alpha v\sb \beta\otimes c\sp{\alpha\beta}  =  u\sb \alpha v\sw
1\otimes \pi(w v\sw 2)\sp\alpha
=u\sw 1v\sw 1\otimes \pi(wu\sw 2v\sw 2) = (uv)\sb \alpha \otimes c\sp \alpha,
$$
where $w\sp \alpha\in \pi\sp{-1}(c\sb \alpha)$ and $w\in\pi\sp{-1}(c)$.
Hence the condition (\ref{diag.A}) is satisfied. Furthermore, we have
$\epsilon(c\sp \alpha)u\sb \alpha = \epsilon(\pi(wu\sw 2))u\sw 1 =
\epsilon(c)u$
and \begin{eqnarray*}
u\sb{\alpha\beta}\otimes c\sw 1\sp \beta\otimes c\sw 2\sp \alpha  & = &
 u\sw 1\sb \alpha\otimes \pi(w\sw 1 u\sw 2)\otimes c\sw 2 \sp \alpha = u\sw
1\otimes
\pi(w\sw 1u\sw 2)\otimes \pi(w\sw 2u\sw 3) \\
& = & u\sw 1\otimes \pi(wu\sw 2)\sw 1\otimes \pi(wu\sw 2)\sw 2 =
u\sb \alpha\otimes c\sp \alpha \sw 1\otimes c\sp \alpha\sw 2,
\end{eqnarray*}
where again $w\sw 1\in\pi\sp{-1}(c\sw 1)$ and $w\sw 2\in\pi\sp{-1}(c\sw 2)$.
Therefore the condition (\ref{diag.B}) is also
satisfied. \endproof

We note that diagrams (\ref{diag.A}) and (\ref{diag.B}) are dual to
each other  in the
following sense. The diagrams (\ref{diag.B}) may be obtained from diagrams
(\ref{diag.A}) by
interchanging $\cdot$ with $\Delta$, $\eta$ with $\epsilon$ and $P$ with
$C$, and by
reversing the arrows. With respect to this duality property the axioms for
the map
$\psi$ are self-dual. Therefore we  can dualise
Proposition~\ref{bundle.prop} to obtain the following:

\begin{prop}
Let $C,P$ be entwined by $\psi:C\tens P\to P\tens C$.  For every algebra
character $\kappa : P\to k$  we have the following:

1. For any positive integer $n$, $C\sp{\otimes n}$ is a right
$P$-module with the action $\ra \sp n = (\kappa\otimes \id\sp n)\circ\psi\sb
{12}\circ\psi\sb
{23}\circ\ldots\circ\psi\sb{ nn+1} = (\kappa\otimes \id\sp
n)\circ\overrightarrow{\psi}\sp n$.

2. The action $\ra\sp n$ maps $\Delta\sp n (C)$ to itself.

3. The  subspace $I_\kappa=\span \{c\ra\sp 1 u-c\kappa(u)|c\in C,\
u\in P\}$ is
a coideal. Hence $M=C/I_\kappa$ is a coalgebra. We denote the canonical
projection by $\pi_\kappa:C\to M$.

4. There is a map $\zeta^M:C\tens P\to C\tens^M C$ defined by $\zeta^M(c\tens
u)=c\o\tens^M c\t\ra\sp 1 u$, where $C\tens^M C=\span\{ c\tens d\in C\tens
C|c\o\tens \pi_\kappa(c\t) \tens d=c\tens \pi_\kappa(d\o)\tens d\t\}$ is the
cotensor product under  $M$. If $\zeta^M$ is a bijection, we say that
$C(M,P,\psi,\kappa)$ is a dual $\psi$-principal bundle.
\label{cobundle.prop}
\end{prop}
\proof 1. The explicit action is
\[ (c_n\tens\cdots\tens c_1)\ra^n u=c_n^{\alpha_n}\tens \cdots\tens
c_1^{\alpha_1} \kappa(u_{\alpha_1\cdots \alpha_n}).\]
Then clearly
\begin{eqnarray*}
((c_n\tens\cdots \tens c_1)\ra^n u)\ra^n v &=&
c_n^{\alpha_n\beta_n}\tens\cdots\tens
c_1^{\alpha_1\beta_1}\kappa(u_{\alpha_1\cdots\alpha_n}
v_{\beta_1\cdots\beta_n})\\
&=& c_n^{\alpha_n}\tens c_{n-1}^{\alpha_{n-1}\beta_{n-1}}\tens\cdots\tens
c_1^{\alpha_1\beta_1}\kappa((u_{\alpha_1\cdots\alpha_{n-1}}
v_{\beta_1\cdots\beta_{n-1}})_{\alpha_n})\\
&=& \cdots = c_n^{\alpha_n}\tens \cdots\tens c_1^{\alpha_1}
\kappa((uv)_{\alpha_1\cdots \alpha_n})=(c_n\tens\cdots\tens c_1)\ra^n(uv)
\end{eqnarray*}
for all $c_i\in C$ and $u,v\in P$. We used (\ref{ind.A}) repeatedly.

2. We have $(c\o\tens c\t)\ra^2 u=c\o{}^\beta\tens
c\t{}^\alpha\kappa(u_{\alpha\beta})=c^\alpha\o\tens
c^\alpha\t\kappa(u_\alpha)=\Delta
(c\ra^1 u)$ by (\ref{ind.B}). Similarly for higher $\Delta^n(C)$.

3. Explicitly,   $I_\kappa=\span\{ c^\alpha\kappa(u_\alpha)-c\kappa(u)| c\in
C,\ u\in P\}$. But using (\ref{ind.B}) we have
$\Delta(c^\alpha\kappa(u_\alpha)-c\kappa(u))=c\o{}^\beta\tens c\t{}^\alpha
\kappa(u_{\alpha\beta})-c\o\tens c\t\kappa(u)=c\o\tens
(c\t{}^\alpha\kappa(u_\alpha)-c\t
\kappa(u))+(c\o{}^\beta\kappa(u_{\alpha\beta})-c\o\kappa(u_\alpha))\tens
c\t{}^\alpha\in C\tens I_\kappa+I_\kappa\tens C$. Hence $I_\kappa$ is a
coideal.

4. The stated map $\zeta^M(c\tens u)=c\o\tens
c\t{}^\alpha\kappa(u_\alpha)$ has its image in $C\tens^M C$ since
\[ c\o\tens \pi_\kappa(c\t )\tens c\th{}^\alpha \kappa(u_\alpha)=c\o\tens
\pi_\kappa(c\t{}^\beta) \kappa(u_{\alpha\beta})\tens c\th{}^\alpha\]
using (\ref{ind.B}) and $\pi_\kappa(I_\kappa)=0$. By dimensions in the
finite-dimensional case, it is natural to require that this is an
isomorphism.   \endproof

This is also an example of a dual version of the theory of $C$-Galois
extensions. The proposition is dual to Proposition~\ref{bundle.prop} in the
sense that all arrows are reversed. In concrete terms, if $P,C$ are
finite-dimensional then
$\psi^*:P^*\tens C^*\to C^*\tens P^*$ and  $\kappa\in P^*$ make
$C^*(M^*,P^*,\psi^*,\kappa)$ a $\psi^*$-principal bundle. Here $M^*=\{f\in
C^*|(\kappa\tens f)\circ\psi= f\tens\kappa\}$. If $C,P$ are entwined and we
have both $e\in C$ and $\kappa:P\to k$, we can have both a $\psi$-principal
bundle and a dual one at the same time. An obvious example, in the setting of
Example~\ref{quantum.ex}, is $P=C=H$ a Hopf algebra and $\psi(c\tens
u)=u\o\tens cu\t$ by the coproduct. Then Proposition~\ref{bundle.prop} with
$e=1$ gives a quantum principal bundle with $M=k$ and right coaction given by
the coproduct. On the other hand, Proposition~\ref{cobundle.prop} with
$\kappa=\eps$ gives a dual bundle with action by right multiplication.

Finally, we note that there is a close connection with the theory of
factorisation of (augmented) algebras introduced   in
\cite{Ma:phy}\cite{Ma:rem}  as part of a factorisation theory of Hopf
algebras.
According to this theory, a factorisation of an algebra $X$ into subalgebras
$A,B$ (so that the product $A\tens B\to X$ is a linear isomorphism) is
equivalent to a {\em factorisation structure} $B\tens A\to A\tens B$ with
certain properties. It was also shown that when $A,B$ are augmented by
algebra
characters then the factorisation structure induces a right action of $A$ on
$B$ and a left action of $B$ on $A$, respectively.

\begin{prop} Let $C$ be finite-dimensional. Then an entwining structure
$\psi:C\tens P\to P\tens C$ is equivalent by partial dualisation to a
factorisation structure $\tilde\psi:C^*\tens P^{\rm op}\to P^{\rm op}\tens
C^*$.  In the augmented case, the induced coaction $\Delta^1_R$ and action
$\ra^1$ in Propositions~\ref{bundle.prop} and~\ref{cobundle.prop} are the
dualisations of the actions induced by the factorisation.
\label{factor.prop}
\end{prop}
\proof We use the notation $\tilde\psi(f\tens u)=u_i\tens f^i$ say, for $f\in
C^*$ and $u\in P$. The equivalence with $\psi$ is by
$u_i\<f^i,c\>=u_\alpha\<f,c^\alpha\>$, where $\<\ ,\ \>$ denotes the
evaluation
pairing. It is easy to see that $\psi$ entwines $C,P$ iff $\tilde\psi$ obeys
\cite{Ma:rem}cf\cite{Ma:phy}
\eqn{factor.A}{
\tilde\psi\circ(\cdot\tens\id)=(\id\tens\cdot)\circ\tilde\psi_{12}
\circ\tilde\psi_{23},\quad \tilde\psi(f\tens 1)=1\tens f}
\eqn{factor.B}{
\tilde\psi\circ(\id\tens\cdot)=(\cdot\tens\id)\circ\tilde\psi_{23}
\circ\tilde\psi_{12},\quad \tilde\psi(1\tens u)=u\tens 1}
for all $f\in C^*$ and $u\in P^{\rm op}$. Thus, the first of these is
$u_i\<c,(fg)^i\>=u_\alpha\<c^\alpha,fg\>=u_\alpha\<c^\alpha\o,f\>
\<c^\alpha\t,g\>=u_{\alpha\beta}\<c\o{}^\beta,f\>\<c\t{}^\alpha,g\>=u_{\alpha
i}\<c\o,f^i\>\<c\t{}^\alpha,g\>=u_{ji}\<c\o,f^i\>\<c\t,g^j\>$ using
(\ref{ind.B}). Similarly for  (\ref{factor.B}) using (\ref{ind.A}), provided
we remember to use the opposite product on $P$. Such data $\psi^{\rm op}$ is
equivalent by \cite{Ma:rem}\cite{Ma:phy}  to the existence of an algebra $X$
factorising into $P^{\rm op}C^*$. Given such $X$ we recover $\tilde\psi$ by
$uc=\cdot\circ \tilde\psi(c\tens u)$ in $X$, and conversely, given
$\tilde\psi$
we define $X=P^{\rm op}\und\tens C^*$ as in (\ref{btens}), but with
$\tilde\psi$.
Also from this theory, if we have $\kappa$ an algebra character on
$P^{\rm op}$
(or on $P$) then $\ra=(\kappa\tens\id)\circ\tilde\psi$ is a right action of
$P^{\rm op}$ on $C^*$, which clearly dualises to the right action of
$P$ on $C$
in Proposition~\ref{cobundle.prop}. Similarly, if $e$ is a character
on $C^*$
then $\la=(\id\tens e)\circ\tilde\psi$ is a left action of $C^*$ on
$P^{\rm op}$
(or on $P$) which clearly dualises to the right coaction of $C$ in
Proposition~\ref{bundle.prop}. \endproof

An obvious setting in which factorisations arise is the braided tensor
product
(\ref{btens}) of algebras in braided
categories\cite{Ma:bra}\cite{Ma:bg}\cite{Ma:introp}, with
$\tilde\psi=\Psi$ the
braiding. Thus if $A\und\tens B$ is a braided tensor product of
algebras (e.g.
of module algebras under a background quantum group) we can look for a
suitable
dual coalgebra $B^*$ in the category and the corresponding entwining
$\psi$ of
$B^*,A^{\rm op}$. This provides a large class of entwining structures.

Another source is the theory of double cross products $G\bowtie H$ of Hopf
algebras in \cite{Ma:phy}. These factorise as Hopf algebras and hence, in
particular, as algebras. In this context,
Proposition~\ref{factor.prop} can be
combined with the result in \cite[Sec.~3.2]{Ma:phy} that the double cross
product is equivalent by partial dualisation to a {\em bicrossproduct}
$H^*\bicross G$. These bicrossproduct Hopf algebras (also due to the second
author) provided one of the first general constructions for
non-commutative and
non-cocommutative Hopf algebras, and many examples are known.

\begin{prop} Let $C\bicross P^{\rm op}$ be a bicrossproduct
bialgebra\cite[Sec.~3.1]{Ma:phy}, where $P^{\rm op}, C$ are bialgebras
suitably
(co)acting on each other. Then $C,P$ are entwined by
\[ \psi(c\tens u)=u\o\sco 0\tens u\o\sco 1 (u\t\la c).\]
Here  $\la$ is the left action of $P^{\rm op}$ on $C$ and $u\sco 0\tens u\sco
1$ is the right coaction of $C$ on $P^{\rm op}$, as part of the
bicrossproduct
construction.
\end{prop}
\proof We derive this result under the temporary assumption that $C$ is
finite-dimensional. Thus the bicrossproduct is equivalent to a double cross
product $P^{\rm op}\bowtie C^*$ with actions $\la,\ra$ defined by $f\la u=u\sco
0\<f,u\sco 1\>$ and $\<u\la c,f\>=\<c,f\ra u\>$ for all $f\in C^*$. Then
$\tilde\psi$ for this factorisation is $\tilde\psi(f\tens u)=f\o\la u\o\tens
f\t\ra u\t$ according to \cite{Ma:phy}\cite{Ma:rem}. The correspondence in
Proposition~\ref{factor.prop} then gives $\psi$ as stated. Once the formula for
$\psi$ is known, one may verify
directly that it entwines $C,P$ given the compatibility conditions between the
action and coaction of a bicrossproduct in \cite[Sec.~3.1]{Ma:phy}. \endproof

Now we describe  trivial $\psi$-principal
bundles and gauge transformations in them.

\begin{prop}
Let $P$ and  $C$ be entwined by $\psi$ as in
Definition~\ref{entwined.df} and let $e$ be a group-like element in
$C$. Assume the following data:

1. A map $\psi\sp C:C\otimes C\to C\otimes C$ such that
\begin{equation}
{(\id\tens\Delta)\circ\psi\sp C=\psi\sp C_{12}\circ\psi^C_{23}
\circ(\Delta\tens\id),\quad (\id\tens\eps)\circ\psi^C=\eps\tens\id,}
\label{psiC.condition}
\end{equation}
and $\psi^C(e\otimes c) = \Delta c$, for any $c\in C$;

2. A convolution invertible map $\Phi: C\to P$ such that $\Phi(e) =1$
and
\begin{equation}
\psi\circ(\id\otimes\Phi) = (\Phi\otimes\id)\circ\psi^C.
\label{cov.phi}
\end{equation}

Then there is a $\psi$-principal bundle over $M=P_e^{coC}$ with
structure coalgebra $C$ and total space $P$.  We call it the
{\em trivial $\psi$-principal bundle}
$P(M,C,\Phi,\psi,\psi^C, e)$  associated to our data, with trivialisation
$\Phi$.
\label{trivial.prop}
\end{prop}

\proof  The proof of the
proposition is similar to the proof that the
trivial quantum  principal bundle  in \cite[Example 4.2]{BrzMa:gau} is in fact
a quantum principal bundle. First we observe that the map
$$
\Theta: M\otimes C\to P, \qquad x\otimes c\mapsto x\Phi(c)
$$
is an isomorphism of linear spaces. Explicitly the inverse is given by
$$
\Theta^{-1} :u\mapsto u\sco 0\Phi\sp{-1}({u\sco 1}\sw 1)\otimes {u\sco 1}\sw 2.
$$
To see that the image of the above map is in $M\otimes C$ we first
notice that (\ref{cov.phi}) implies that $\Delta^1\sb R \circ\Phi =
(\Phi\otimes\id)\circ\Delta$ and that
\begin{equation}
\phi(c\sw 1\otimes
\Phi^{-1}(c\sw 2)) = \Phi^{-1}(c)\otimes e .
\label{cov.phi-1}
\end{equation}
Therefore for any $u\in P$,
$$
\Delta^1\sb R(u\sco 0\Phi^{-1}(u\sco 1)) = u\sco 0\psi({u\sco 1}\sw 1
\otimes \Phi^{-1}({u\sco 1}\sw 2)) = u\sco 0\Phi^{-1}(u\sco 1)\otimes e,
$$
and thus $ u\sco 0\Phi^{-1}(u\sco 1)\in M$. Then it is easy to prove
that the above maps are inverses to each other.

We remark that $\Theta$ is in fact a left $M$-module and a right
$C$-comodule map, where the coaction in $M\otimes C$ is given by
$x\otimes c \mapsto x\otimes c\sw 1\otimes c\sw 2$. Moreover
$\psi\circ(\id\otimes \Theta) = (\Theta\otimes\id)\circ
\psi^C_{23}\circ\psi_{12}$.

The proof that $\chi\sb M$ in this case is a bijection follows exactly
the method used in the proof of \cite[Example~4.2]{BrzMa:gau} and thus
we do not repeat it here. \endproof\vspace{12pt}

Next, we consider gauge transformations.

\begin{df} Let $P(M,C,\Phi,\psi,\psi^C, e)$ be a trivial
$\psi$-principal bundle as in Proposition~\ref{trivial.prop}. We say
that a convolution invertible map $\gamma :C\to M$ such that
$\gamma(e) =1$ is a {\em gauge transformation} if
\begin{equation}
\psi^C_{23}\circ\psi_{12}\circ(\id\otimes\gamma\otimes\id)\circ(\id\tens\Delta)
=
(\gamma\otimes\id\otimes\id)\circ(\Delta\otimes\id)\circ\psi^C .
\label{gauge.condition}
\end{equation}
\label{gauge.df}
\end{df}
\begin{prop}
If $\gamma : C\to M$ is a gauge transformation in
$P(M,C,\Phi,\psi,\psi^C, e)$ then $\Phi' = \gamma*\Phi$, where $*$
denotes the convolution product,  is a
trivialisation of $P(M,C,\Phi,\psi,\psi^C, e)$. The set of all gauge
transformations in $P(M,C,\Phi,\psi,\psi^C, e)$ is a group with
respect to the convolution product. We say that two trivialisations
$\Phi$ and $\Phi'$ are {\em gauge equivalent} if there exists a gauge
transformation $\gamma$ such that $\Phi' = \gamma *\Phi$.
\label{gauge.prop}
\end{prop}
\proof Clearly $\Phi'$ is a convolution invertible map such that
$\Phi'(e) =1$. To prove that it satisfies (\ref{cov.phi}) we first
introduce the notation
$$
\psi^C(b\otimes c) = c_A\otimes b^A \qquad ({\rm summation \;\;\;assumed}),
$$
in which the condition (\ref{gauge.condition}) reads explicitly
$$
\gamma(c\sw 1)_\alpha \otimes {c\sw 2}_A \otimes b^{\alpha A} =
\gamma({c_A}\sw 1)\otimes {c_A}\sw 2 \otimes b^A,
$$
and then compute
\begin{eqnarray*}
\psi(\id\otimes\Phi')(b\otimes c) & = & \psi(b\otimes \gamma(c\sw
1)\Phi(c\sw 2))   =  \gamma(c\sw 1)_\alpha
\Phi(c\sw 2)_\beta\otimes b^{\alpha\beta} \\
& = & \gamma(c\sw 1)_\alpha \Phi({c\sw 2}_A)\otimes b^{\alpha A}
=  \gamma({c_A}\sw 1)\Phi( {c_A}\sw 2) \otimes b^A \\
& = & (\Phi'\otimes\id)\psi^C(b\otimes c).
\end{eqnarray*}
This proves the first part of the proposition.

Assume now that $\gamma_1$, $\gamma_2$ are gauge transformations. Then
\begin{eqnarray*}
(\gamma_1(c\sw 1)\gamma_2(c\sw 2))_{\alpha}\otimes c\sw 3_A \otimes
b^{\alpha A} & = & \gamma_1(c\sw 1)_\alpha \gamma_2(c\sw 2)_{\beta}\otimes c\sw
3_A \otimes
b^{\alpha\beta A} \\
& = &\gamma_1(c\sw 1)_\alpha \gamma_2(c\sw 2_A\sw 1)\otimes c\sw
2_A\sw 2 \otimes b^{\alpha A}\\
& = & \gamma_1(c_A\sw 1) \gamma_2(c_A\sw 2)\otimes c_A\sw 3
\otimes b^ A.
\end{eqnarray*}
Therefore $\gamma_1*\gamma_2$ is a gauge transformation too. Clearly
$\epsilon$ is a gauge transformation and thus provides the
unit. Finally, to prove that if $\gamma$ is a gauge transformation
then so is $\gamma^{-1}$, we observe that if $\gamma_3 =
\gamma_1*\gamma_2$ and $\gamma_2$ are gauge transformations then so is
$\gamma_1$. Indeed, if $\gamma_1*\gamma_2$ is a gauge transformation
then
$$
(\gamma_1(c\sw 1)\gamma_2(c\sw 2))_{\alpha}\otimes c\sw 3_A \otimes
b^{\alpha A} = \gamma_1(c_A\sw 1)\gamma_2(c_A\sw 2)\otimes c_A\sw 3
\otimes b^ A ,
$$
but since $\gamma_2$ is a gauge transformation, we obtain
$$
\gamma_1(c\sw 1)_\alpha \gamma_2(c\sw 2_A\sw 1)\otimes c\sw
2_A\sw 2 \otimes b^{\alpha A}
=  \gamma_1(c_A\sw 1) \gamma_2(c_A\sw 2)\otimes c_A\sw 3
\otimes b^ A.
$$
Applying $\gamma_2^{-1}$ to the second factor in the tensor product
and then multiplying the first two factors we obtain
$$
\gamma_1(c\sw 1)_\alpha \otimes {c\sw 2}_A \otimes b^{\alpha A} =
\gamma_1({c_A}\sw 1)\otimes {c_A}\sw 2 \otimes b^A,
$$
i.e. $\gamma_1$ is a gauge transformation as stated. Now applying this
result to $\gamma_3 = \epsilon$ and $\gamma_2 =\gamma$ we deduce that
$\gamma^{-1}$ is a gauge transformation as required. This completes
the proof of the
proposition. \endproof \vspace{12pt}

Although the existence of the map $\psi^C$ as in
Proposition~\ref{trivial.prop} is not guaranteed for all
coalgebras, the map $\psi^C$ exists in all the examples discussed in
this section:

\begin{ex}\label{quantum.trivial} For a quantum principal bundle $P(M,H)$ as in
Example~\ref{quantum.ex}, we define
\[\psi^H(b\otimes c)=c\sw 1\otimes bc\sw 2, \]
for all $b,c\in H$. Then (\ref{trivial.prop})--(\ref{gauge.prop}) reduces to
the theory of trivial quantum principal bundles and their gauge transformations
in \cite{BrzMa:gau}.
\end{ex}
\proof It is easy to see by standard Hopf algebra calculations that
(\ref{psiC.condition}) is satisfied by the bialgebra axiom for $H=C$ in this
case. Moreover, (\ref{cov.phi}) reduces to $\Phi$ being an intertwiner of
$\Delta_R$ with $\Delta$. The condition (\ref{gauge.condition}) is empty. This
recovers the setting introduced in \cite{BrzMa:gau}. \endproof\vspace{12pt}

In the braided case we use the above theory to  arrive at a natural definition
of trivial braided principal bundle:

\begin{ex}\label{braided.trivial} For a braided
principal bundle $P(M,B,\Psi)$ as in Example~\ref{braided.ex}, we define
a trivialisation as a convolution-invertible unital morphism $\Phi:B\to P$ in
the braided category such that
$\Delta_R\circ\Phi=(\Phi\tens\id)\circ \und\Delta$, where $\Delta_R$ is the
braided right coaction of $B$ on $P$. We define a gauge transformation
as a convolution-invertible unital morphism $\gamma:B\to M$, acting on
trivialisations by the convolution product $*$. This is a trivial
$\psi$-principal bundle with
\[ \psi^B(b\tens c)=\Psi(b\otimes c\Bo)c\Bt,\]
where $\und\Delta c=c\Bo\tens c\Bt$ is the braided group coproduct.
\end{ex}
\proof This time, (\ref{psiC.condition}) follows from the braided-coproduct
homomorphism property of a braided group\cite{Ma:bra}. From this and the form
of $\psi$, we see that (\ref{cov.phi}) becomes
$\Delta_R\circ\Phi(c)=(\Phi\tens\id)\circ\Psi(b\tens c\Bo)c\Bt$. Setting $b=e$
gives the condition stated on $\Phi$ because the braiding with $e=1$ is always
trivial. Assuming the stated condition, (\ref{cov.phi}) then becomes
$\Phi(c\Bo)\tens b c\Bt=(\Phi\tens \id)\circ\Psi(b\tens c\Bo)c\Bt$,
which is equivalent (by replacing $c\Bt$ by $c\Bt\tens \und S c\Bth$ and
multiplying, where $\und S$ is the braided antipode) to
\[ (\Phi\otimes\id)\circ\Psi = \Psi\circ(\id\otimes\Phi).\]
When all our constructions take place in a braided category,  this is the
functoriality property  implied by requiring that $\Phi$ is a morphism in the
category. The theory of trivial $\psi$-bundles only requires this functoriality
condition itself. Similarly, we compute the gauge condition
(\ref{gauge.condition}) using  $\psi(b\tens \gamma(c))=\Psi(b\tens\gamma(c))$
because $\gamma(c)\in M$, and operate on it by replacing $c\Bt$ by $c\Bt\tens
\und S c\Bth$ and multiplying. Then it reduces to
$\Psi_{23}\circ\Psi_{12}\circ(\id\tens
\gamma\circ\und\Delta)=(\gamma\circ\und\Delta\tens\id)\circ\Psi$.
Since $\und\Delta$ is a morphism, we see (by applying the braided group counit)
that the gauge condition (\ref{gauge.condition}) is equivalent to
\[ (\gamma\otimes\id)\circ\Psi = \Psi\circ(\id\otimes\gamma).\]
As before, this is naturally implied by requiring that $\gamma$ is a morphism
in our braided category. It is clear that the convolution product $*$ preserves
the property of being a morphism since $\und\Delta$ and  $\Delta_R$ are assumed
to be morphisms.  \endproof\vspace{12pt}

Also, for a $\psi$-principal
bundle over a quantum homogeneous space as in
Example~\ref{embeddable.ex},  we define
\begin{equation}\label{psiC.embed} \psi^C(b\tens c)= \pi(v\sw 1)\otimes\pi(u
v\sw 2),\end{equation}
where $u\in\pi^{-1}(b)$, $v\in\pi^{-1}(c)$. Then a trivialisation of the
bundle is a convolution-invertible map $\Phi:C\to H$ obeying
$\Phi\circ\pi(1)=1$ and
\begin{equation}\label{Phi.embed} \Phi(c)\o\tens
\pi(u\Phi(c)\t)=\Phi\circ\pi(v\o)\tens \pi(uv\t)\end{equation}
for all $c\in C,u\in H$, and $v\in \pi^{-1}(c)$. Taking $u=1$ requires, in
particular, the natural intertwiner condition
$(\Phi\tens\id)\circ\Delta=\Delta_R\circ \Phi$. There is, similarly, a
condition on gauge transformations $\gamma$ obtained from
(\ref{gauge.condition}). Hence our formulation of trivial $\psi$-principal
bundles covers all the main sources of $\psi$-principal bundles discussed in
this section.

We conclude this section with some explicit examples of $\psi$-principal
bundles.

\begin{ex}\rm   Let $H$ be
a quantum cylinder $A\sb q\sp{2|0}
[x\sp{-1}]$, i.e. a free associative algebra generated by $x, x\sp{-1}$
and $y$ subject to the relations $yx=qxy$, $xx\sp{-1} = x\sp{-1}x =1$,
with a natural Hopf algebra structure:
\begin{equation}
\Delta x\sp{\pm 1} = x\sp{\pm 1}\otimes x\sp{\pm 1}, \qquad
\Delta y = 1\otimes y +y\otimes x, \quad etc.
\label{coproduct.qplane}
\end{equation}
Consider a right ideal $J$ in $H$ generated by  $x-1$ and $x\sp{-1}-1$.
Clearly, $J$ is a coideal and therefore $C=A\sb q\sp{2|0}[x\sp{-1}]/J$
is a coalgebra and a canonical epimorphism $\pi :H\to C$ is a coalgebra
map. $C$ is spanned by the elements $c\sb n = \pi(y\sp n)$,
$n\in{\bf Z}\sb{\geq 0}$, and the coproduct and the counit are given
by
\begin{equation}
\Delta{c\sb n} = \sum\sb{k=0}\sp{n}\qbinom nkq
c\sb k\otimes c\sb{n-k},\qquad \eps(c_n)=0.
\label{coproduct.braided.line}
\end{equation}
 We are in the situation of
Example~\ref{embeddable.ex} and thus
we have the entwining structure $\psi : C\otimes H\to H\otimes C$, which
explicitly  computed comes out as
\begin{equation}
\psi(c\sb l\otimes x\sp my\sp n) = \sum\sb{k=0}\sp{n}\qbinom nkq
q\sp{l(k+m)}x\sp m y\sp k\otimes c\sb{n+l-k},
\label{cylinder.psi.eqn}
\end{equation}
where
$$
\sum\sb{k=0}\sp{n}\qbinom nkq = {\frac{[n]\sb q!}
{[n-k]\sb q! [k]\sb q!}},
$$
and
$$
[n]\sb q! = [n]\sb q\cdots [2]\sb
q[1]\sb q, \quad [0]\sb q! =1,
\qquad [n]\sb q = 1+q+\ldots+q\sp{n-1}.
$$
{}From this definition of $\psi$ one easily computes the right coaction
of $C$ on $H$ as well as the fixed point subalgebra $M=k[x,x\sp{-1}]$,
i.e. the algebra of functions on a circle.
By Example~\ref{embeddable.ex} we have just constructed a generalised
quantum principal bundle $A\sb q\sp{2|0}[x\sp{-1}](k[x,x\sp{-1}],C,
\psi,c\sb 0)$.

Finally we note that the above bundle is trivial in the sense of
Proposition~\ref{trivial.prop}. The trivialisation
$\Phi:C\to A\sb q\sp{2|0}[x\sp{-1}]$
and its inverse $\Phi\sp{-1}$ are defined by
\begin{equation}
\Phi(c\sb n) = y\sp n ,\qquad \Phi\sp{-1}(c\sb n) = (-1)\sp
nq\sp{n(n-1)/2}y\sp n .
\label{trivialisation.eqn}
\end{equation}
One can easily check that the map $\Phi$ satisfies the required
conditions. Explicitly, the map $\psi^C :C\otimes C\to C\otimes C$
reads
$$
\psi^C(c_m\otimes c_n) = \sum\sb{k=0}\sp{n}\qbinom nkq q^{km}
c\sb k\otimes c\sb{m+n-k}.
$$
Therefore
$$
\psi\circ(\id\otimes\Phi)(c_m\otimes c_n) = \psi(c_m\otimes y^n) =
\sum\sb{k=0}\sp{n}\qbinom nkq q^{km}
y\sb k\otimes c\sb{m+n-k} = (\Phi\otimes\id)\circ\psi^C(c_m\otimes
c_n).
$$
Since the bundle discussed in this example is trivial, we can compute
its gauge group. One easily finds that a convolution invertible map
$\gamma:C\to k[x,x^{-1}]$ satisfies condition (\ref{gauge.condition})
if and only if
$\gamma(c_n) = \Gamma_n x^n$ (no summation), where $n\in {\bf Z}_{\geq 0}$,
$\Gamma_n\in k$ and
$\Gamma_0 =1$. Therefore the gauge group is equivalent to the group of
sequences $\Gamma = (1, \Gamma_1,
\Gamma_2, ...)$ with the product given by
$$
(\Gamma\cdot\Gamma')_n = \sum\sb{k=0}\sp{n}\qbinom nkq
\Gamma\sb k\Gamma'\sb{n-k}.
$$
\label{cylinder.ex}
\end{ex}

For the simplest example of a braided principal bundle, one can simply take any
braided group $B$ and any algebra $M$ in the same braided category. Then the
braided tensor product algebra $P= M\und\tens B$, along with the definitions
\begin{equation} \label{tensor.braided} \Delta_R=\id\tens\und\Delta,\quad
\Phi(b)=1\tens b,\quad \Phi^{-1}(b)=1\tens \und S b\end{equation}
put us in the setting of Examples~\ref{braided.ex} and~\ref{braided.trivial}.
Note first that
$\Delta_R$ is a coaction (the tensor product of the trivial coaction and the
right coregular coaction) and makes $P$ into a braided comodule algebra.
Moreover, the induced map
\[ \chi_M(m\tens b\tens n\tens c)= m\Psi(b\tens n)c\Bo\tens c\Bt\]
for $m,n\in M$, $b,c\in B$, is an isomorphism $P\tens_M P\to P\tens P$; it has
inverse
\[ \chi_M^{-1}(m\tens b\tens c)=m\tens b\und S c\Bo\tens 1\tens c\Bt.\]
It is also clear that $\Phi$ is a trivialisation. This is truly a trivial
braided principal bundle because $P$ is just a (braided) tensor product
algebra.

\begin{ex}\label{line.ex}\rm  Let $B=k[c]$ be the braided line generated by $c$
with braiding $\Psi(c\tens c)=qc\tens c$ and the linear coproduct $\und\Delta
c=c\tens 1+1\tens c$. It lives in the braided category ${\rm Vec}_q$ of
$\Z$-graded vector spaces with braiding $q^{\deg(\ )\deg(\ )}$ times the usual
transposition. Here $\deg(c)=1$. Let $M=k[x,x^{-1}]$ be viewed as a $\Z$-graded
algebra as well, with $\deg(x)=1$. Then $P=k[x,x^{-1}]\und\tens k[c]$ is a
trivial braided principal bundle with the coaction and trivialisation
\begin{equation}\label{line.coaction} \Delta_R( x^m\tens c^n)=\sum_{k=0}^{n}
\qbinom nkq x\sp m \tens c\sp k\otimes c^{n-k},\quad \Phi(c^n)=1\tens
c^n.\end{equation}
As a $\psi$-principal bundle, this example clearly coincides with the preceding
one, albeit constructed quite differently: we identify $c^n=c_n$ and $y=1\tens
c$, and note that in the braided tensor product algebra $k[x,x^{-1}]\und\tens
k[c]$ we have the product $(1\tens c)(x\tens 1)=\Psi(c\tens x)=q(x\tens
1)(1\tens c)$, i.e. $P=A_q^{2|0}[x^{-1}]$. It is also clear that the coproduct
deduced in (\ref{coproduct.braided.line}) can be identified with the braided
line coproduct which is part of our initial data here. This particular braided
tensor product algebra $P$ is actually the algebra part of the bosonisation of
$B=k[c]$ viewed as living in the category of comodules over $k[x,x^{-1}]$ as a
dual-quasitriangular Hopf algebra (see \cite[p510]{Ma:book}), and becomes in
this way a Hopf algebra. This bosonisation is the Hopf algebra $H$ which was
part of the initial data in the preceding example. Finally, gauge
transformations $\gamma$ from the braided point of view are arbitrary
degree-preserving unital maps $k[c]\to k[x,x^{-1}]$, i.e. given by the group of
sequences $\Gamma$ as found before.  \end{ex}

This example demonstrates the strength of braided group gauge theory; even the
most trivial braided quantum principal bundles may be quite complicated when
constructed by more usual Hopf algebraic means. On the other hand, the
following embeddable quantum homogeneous space does not appear to be of the
braided type, nor (as far as we know) a trivial bundle.

\begin{ex}
Let $P$ be the algebra of functions on the quantum group $GL_q(2)$. This
is generated by elements $\alpha$, $\beta$, $\gamma$, $\delta$ and $D$
subject to the relations
$$
\alpha\beta = q\beta\alpha, \qquad \alpha\gamma = q\gamma\alpha, \qquad
\alpha\delta = \delta\alpha + (q-q^{-1})\beta\gamma ,\qquad
\beta\gamma = \gamma\beta ,
$$
$$
 \beta\delta = q\delta\beta , \qquad
\gamma\delta = q\delta\gamma ,\qquad
(\alpha\delta - q\beta\gamma)D = D(\alpha\delta - q\beta\gamma) =
1.
$$
Let $C$ be a vector  space spanned by $c_{m,n}$, $m\in {\bf Z}_{>0}$, $n\in\bf
Z$
with the coalgebra structure
$$
\Delta(c_{i,j}) = \sum_{k=0}^mq^{k(m-k)}\qbinom  mk\q2c_{k,n}\otimes
c_{m-k,n+k}, \qquad \epsilon(c_{m,n}) = \delta_{m0} .
$$
Let the linear map $\psi :C\otimes P\to P\otimes C$ be given by
\begin{eqnarray}
\psi(c_{i,j}\otimes\alpha^k\gamma^l\beta^m\delta^n D^r) & = &
\sum_{s=0}^m\sum_{t=0}^n \qbinom ms\q2 \qbinom nt\q2 q^{(m-s)(s+t-l)
+(n-t)t -i(k+l-t-s)}\times \nonumber\\
&& \alpha^{k+m-s}\gamma^{l+n-t}\beta^s\delta^t D^r \otimes
c_{i+m+n-s-t,j+r+t+s}.\label{psi.glq2}
\end{eqnarray}
Then $\psi$ entwines $P$ with $C$. Furthermore if we take $e=a_{0,0}$
then the fixed
point subalgebra $P_e^C$ is generated by $1,\alpha,\gamma$ and hence it
is isomorphic to $A_{1/q}^{2|0}$ and there
is a $\psi$-principal bundle $P(A^{2|0}_{1/q}, C,\psi, e)$.
\label{glq2.ex}
\end{ex}
\proof We reformulate the $C$-Galois extension or $\rho$-principal bundle from
Section~5 of \cite{Brz:hom} as a $\psi$-principal bundle. Using the Part 5 of
Proposition~\ref{bundle.prop} and $\rho$ from \cite{Brz:hom} one can  easily
obtain $\psi$ as stated. In \cite{Brz:hom} it is also noted that the coalgebra
$C$ can be equipped with the algebra structure of
$A_{q^{-2}}^{2|0}[x^{-1}]$  by setting $c_{m,n} = q^{-mn}x^ny^m$. The
coproduct in $C$ is then the same as in Example~\ref{cylinder.ex},
Eq.~(\ref{coproduct.qplane}).
\endproof

\section{Connections in the Universal Differential Calculus Case}
{}From the first assertion of
Proposition~\ref{bundle.prop} we know that the natural
coaction $\Delta\sb R = \Delta\sb R\sp 1$ of $C$ on $P$ extends to the
coaction of $C$ on the tensor product algebra $P\sp{\otimes n}$ for
any positive integer
$n$. Still most importantly this coaction can be restricted to
$\Omega\sp n P$ by the second assertion of
Proposition~\ref{bundle.prop}. Therefore the coalgebra
$C$ coacts on the algebra of
universal forms on $P$. The universal differential structure on $P$ is
covariant with respect to the coaction $\Delta\sb R\sp n$ in the
following sense
\begin{prop}
Let $P$, $C$, $\psi$ and $e$ be as in
Proposition~\ref{bundle.prop}. Let $\d:P\to\Omega\sp 1P$
be the Karoubi differential, $\d u = 1\otimes u-u\otimes 1$ extended
to the whole of $\Omega P$ in the standard way. Then $\overleftarrow{\psi}\sp
n\circ(\id\otimes\d) =(\d\otimes \id)\circ\overleftarrow{\psi}\sp{n-1}$
for any integer $n>1$. Therefore $\Delta\sp n\sb R\circ(\id\otimes \d)
=(\d\otimes \id)\circ \Delta\sb R\sp{n-1}$.
\end{prop}
\proof First we prove the proposition for zero- and one-forms on
$P$. Using the conditions (\ref{diag.A}), for any $u\in P$ and any $c\in C$ we
find
\begin{eqnarray*}
 \overleftarrow{\psi}\sp 2 (c\otimes \d u) & = & \overleftarrow{\psi}\sp
2(c\otimes 1\otimes u) - \overleftarrow{\psi}\sp 2(c\otimes u\otimes
1) \\
&= & (1\otimes u\sb \alpha\otimes c\sp \alpha - u\sb \alpha\otimes 1\otimes
c\sp \alpha) = (\d\otimes \id)\circ\psi(c\otimes u).
\end{eqnarray*}

We represent a general element of $\Omega\sp{n-1}P$ as $\upsilon =
u\sp 1\otimes
u\sp 2\otimes\ldots\otimes u\sp n$. Actually, a general element of
$\Omega\sp{n-1}P$ is a linear combination of the $\upsilon$ such that
the multiplication applied to any two consecutive factors gives
zero. We then note that
$$
\d\upsilon = \d u\sp 1\otimes u\sp
2\otimes\ldots\otimes u\sp n = 1\otimes u\sp 1\otimes u\sp
2\otimes\ldots\otimes u\sp n - u\sp 1\otimes 1\otimes u\sp
2\otimes\ldots\otimes u\sp n,
$$
and that $\overleftarrow{\psi}\sp n = (\overleftarrow{\psi}\sp 2\otimes
\id\sp{n-2})\circ(\id\sp 2\otimes\overleftarrow{\psi}\sp{n-2})$, and
perform calculation similar to the one performed in the zero- and
one-form case to conclude the result.
\endproof

To discuss a theory of connections in $P(M,C,\psi,e)$ it is important
that the horizontal one forms $P\Omega\sp 1MP$ be covariant under the
action of $\Delta\sp 2\sb R$ or, more properly,
$\overleftarrow{\psi}\sp 2$. The following lemma gives a criterion for
the covariance of horizontal one-forms.

\begin{lemma}
For a $\psi$-principal bundle $P(M,C,\psi,e)$ assume that
$\psi(C\otimes M) \subset M\otimes C$. Then $
\overleftarrow{\psi}\sp 2(C\otimes P\Omega\sp 1MP)\subset P\Omega\sp
1MP\otimes C$.
 \label{sufficient.lemma}
\end{lemma}
\proof
Using (\ref{ind.A}) one easily finds that for any $u,v\in P$,
$x,y\in M$ and $c\in C$,
$$
\overleftarrow{\psi}\sp 2(c\otimes ux\otimes yv) =
u\sb \alpha x\sb \beta\otimes y\sb \gamma v\sb \delta\otimes
c\sp{\alpha\beta\gamma\delta}.
$$
If we assume further that $x\sb \beta, y\sb \gamma\in M$ then the result
follows.
\endproof

 We will see later that the hypothesis of Lemma~\ref{sufficient.lemma} is
automatically satisfied for braided principal bundles of
Example~\ref{braided.ex}. In contrary, it is not
necessarily satisfied for $\psi$-bundles on quantum embeddable
homogeneous spaces of Example~\ref{embeddable.ex}. For example, one
can easily check that it is
satisfied
for the bundle discussed in Example~\ref{glq2.ex}. On the other hand  the
$\psi$-principal bundle over the quantum hyperboloid, which is an
embeddable homogeneous space of $E_q(2)$ \cite{BonCic:fre} fails to
fulfil requirements of
Lemma~\ref{sufficient.lemma}.

The covariance of $\Omega P$ and  $P\Omega\sp 1 MP$
enables us to define a connection in
$P(M,C,\psi ,e)$ in a way similar to the definition of a connection in
a quantum principal bundle $P(M,H)$ (compare \cite{BrzMa:gau}).
\begin{df}
Let $P(M,C,\psi ,e)$ be a generalised quantum principal bundle such
that $\psi(C\otimes M)\subset M\otimes C$. A connection  in
$P(M,C,\psi ,e)$ is a left $P$-module
projection $\Pi :\Omega\sp 1 P\to \Omega\sp 1 P$
such that $\ker \Pi = P\Omega\sp 1 MP$ and $\overleftarrow{\psi}\sp
2(\id\otimes\Pi)  = (\Pi
\otimes \id)\overleftarrow{\psi}\sp 2$.
\label{connection.df}
\end{df}

It is clear that for a usual quantum principal bundle $P(M,H)$
Definition~\ref{connection.df} coincides with the definition
of a connection given in \cite{BrzMa:gau}. Thus, the condition in
Lemma~\ref{sufficient.lemma} always holds for
$\psi$ as in Example~\ref{quantum.ex}, while
$\overleftarrow{\psi}\sp 2(\id\otimes\Pi)  = (\Pi
\otimes \id)\overleftarrow{\psi}\sp 2$ if and only if $\Delta\sp
2\sb R\Pi  = (\Pi\otimes \id)\Delta\sb R\sp 2$, which was the condition in
\cite{BrzMa:gau}.

In what follows we assume that the condition in
Lemma~\ref{sufficient.lemma} is satisfied.
A connection $\Pi$ in $P(M,C,\psi, e)$ can be equivalently described
as follows.  First
we define a map $\phi : C\otimes P\tens \ker\eps \to
P\otimes\ker\epsilon\otimes C$
by the commutative diagram
$$
\begin{picture}(200,125)(20,1)
 \put(10,100){\makebox(0,0){$C\otimes \Omega\sp 1 P$}}
 \put(45,100){\vector(1,0){98}}
\put(87,105){$\overleftarrow{\psi}\sp 2$}
\put(170,100){\makebox(0,0){$\Omega\sp 1 P\otimes C$}}
 \put(170,90){\vector(0,-1){50}}
 \put(180,65){$\chi\otimes \id$}
 \put(170,35){\makebox(0,0){$P\otimes \ker\epsilon\otimes C$}}
 \put(20,90){\vector(0,-1){50}}
 \put(30,65){$\id\otimes\chi$}
\put(10,35){\makebox(0,0){$C\otimes P\otimes \ker\epsilon$}}
 \put(50,35){\vector(1,0){85}}
\put(100,43){\makebox(0,0){$\phi$}}
\end{picture} $$
where $\chi(u\otimes v) = u\psi(e\otimes v)$.  The map $\phi$ is
clearly well-defined. Indeed, because $\chi\sb M$ is a bijection,
$\ker\chi = P\Omega\sp 1MP$ and then $\overleftarrow{\psi}\sp
2(C\otimes\ker\chi) \subset
\ker\chi\otimes C$, by
Lemma~\ref{sufficient.lemma}. Therefore $\phi (0) = 0$.

By definition of $P(M,C,\psi,e)$ we have a short exact sequence of
left $P$-module maps
\begin{equation}
0\rightarrow P\Omega\sp 1MP\rightarrow \Omega\sp 1P
\stackrel{\chi}{\rightarrow}
P\otimes\ker\epsilon \rightarrow 0
\label{exact.sequence.universal}
\end{equation}
The exactness of the above sequence is clear since the fact that
$\chi\sb M$ is bijective implies that $\chi$ is surjective and
$\ker\chi = P\Omega\sp 1MP$.  By definition,   $\chi$
intertwines
$\overleftarrow{\psi}\sp 2$ with $\phi$.
\begin{prop}
The existence of a connection $\Pi$ in $P(M,C,\psi ,e)$ is equivalent
to the existence of a left $P$-module splitting $\sigma :P\otimes
\ker\epsilon \to \Omega\sp 1 P$ of the above sequence such that
$\overleftarrow{\psi}\sp 2\circ(\id\otimes\sigma) = (\sigma\otimes
\id)\circ \phi$.
\label{proposition.splitting.universal}
\end{prop}
\proof Clearly the existence of a left $P$-module projection is
equivalent to the existence of a left $P$-module splitting. It remains
to check the required covariance properties. Assume that $\sigma$ has
the required properties then
$$
\overleftarrow{\psi}\sp 2(\id\otimes\Pi) = \overleftarrow{\psi}\sp
2\circ(\id\otimes\sigma)\circ(\id\otimes\chi) = (\sigma\otimes
\id)\circ\phi\circ(\id\otimes\chi) = (\sigma\circ\chi\otimes \id)\circ
\overleftarrow{\psi}\sp 2 =(\Pi\otimes \id)\circ\overleftarrow{\psi}\sp
2.
$$
Conversely, if $\Pi$ has the required properties then one easily finds
that
$$
\overleftarrow{\psi}\sp 2\circ(\id\otimes\sigma\circ\chi) =
(\sigma\otimes \id)\circ\phi\circ(\id\otimes\chi) .
$$
Since $\chi$ is a projection the required property of $\sigma$ follows.
\endproof

To each connection we can associate its connection one-form $\omega :
\ker\epsilon\to \Omega\sp 1  P$ by setting $\omega(c) = \sigma(1\otimes c)$.
\footnote{We can also think of a connection 1-form as a map $H\to
\Omega\sp 1 P$
given by $\omega(c-\epsilon(c))$. This was the  point of
view adopted
in \cite{BrzMa:gau}.} Similarly to the quantum bundle case of
\cite{BrzMa:gau} we have
\begin{prop}
Let $\Pi$ be a connection on $P(M,C,\psi,e)$. Then, for all
$c\in\ker\epsilon$, the  connection
1-form $\omega :\ker\epsilon\to\Omega\sp 1 P$ has the following
properties:

1. $\chi\circ\omega(c) = 1\otimes c$,

2. For any $b,c\in C$, $\overleftarrow{\psi}\sp 2(b\otimes\omega(c)) =
c\su 1{}\sb \alpha c\su 2\sb {\beta\gamma}  \omega(e\sp \gamma)\otimes
b\sp{\alpha\beta}, $
where $c \su 1\tens\sb M c\su 2$ (summation understood) denotes the
{\em translation
map} $\tau(c)=\chi\sb M\sp{-1}(1\otimes
c)$ in $P(M,C,\psi,e)$.

Conversely, if $\omega$ is any linear map $\omega :\ker\epsilon\to
\Omega\sp 1P$ obeying conditions 1-2, then there is a unique
connection $\Pi = \cdot\circ(\id\otimes\omega)\circ\chi$ in $P(M,C,\psi, e)$
such
that $\omega$ is its connection 1-form.
\label{connection.form.prop}
\end{prop}
\proof
For any $b\otimes u\otimes c\in C\otimes P\otimes \ker\epsilon$ the
map $\phi$ is explicitly given by
$$
\phi(b\otimes u\otimes c) = u\sb \alpha c\su 1\sb
\beta c\su 2\sb{\gamma\delta}\otimes e\sp \delta\otimes
b\sp{\alpha\beta\gamma}.
$$
Therefore if $\omega$ is a connection one-form then
\begin{eqnarray*}
\overleftarrow{\psi}\sp 2(b\otimes\omega(c))
&=&\overleftarrow{\psi}\sp 2\circ(\id\otimes\sigma)(b\otimes 1\otimes
c)\\
& = & \sigma(c\su 1\sb \alpha c\su 2\sb{\beta\gamma}\otimes e\sp
\gamma)\otimes b\sp{\alpha\beta} \\
& = &c\su 1\sb
\alpha c\su 2\sb{\beta\gamma}  \omega(e\sp \gamma)\otimes b\sb{\alpha\beta}.
\end{eqnarray*}

Conversely, if $\omega :\ker\epsilon\to\Omega\sp 1P$ satisfies
condition 1 then $\sigma = (\cdot\otimes \id)\circ (\id\otimes\omega)$
gives a left $P$-module splitting of
(\ref{exact.sequence.universal}). Furthermore, Condition~2 implies
\begin{eqnarray*}
(\sigma\otimes \id)\circ\phi(b\otimes u\otimes c) & =
& \sigma(u\sb \alpha c\su 1\sb
\beta c\su 2\sb{\gamma\delta}\otimes e\sp \delta)\otimes
c\sp{\alpha\beta\gamma}\\
& = &   u\sb \alpha c\su 1\sb
\beta c\su 2\sb{\gamma\delta}\omega(e\sp \delta)\otimes
c\sp{\alpha\beta\gamma}
\\
& = &  u\sb \alpha \overleftarrow{\psi}\sp 2(b\sp
\alpha\otimes\omega(c)) = (\cdot\otimes \id)\circ\overleftarrow{\psi}\sp
3(b\otimes u\otimes\omega(c))\\
& = & \overleftarrow{\psi}\sp
2(b\otimes u\omega(c)) = \overleftarrow{\psi}\sp
2\circ(\id\otimes\sigma)(b\otimes u\otimes c).
\end{eqnarray*}
\endproof
\begin{ex}
For a quantum principal bundle $P(M,H)$, Condition~2 in
Proposition~\ref{connection.form.prop} is equivalent to the
$Ad\sb R$-covariance of $\omega$.
\end{ex}
\proof
Using the definition of $\psi$ in Example~\ref{quantum.ex}
one finds
\begin{eqnarray*}
  c\su 1\sb \alpha{c\su 2}\sb{\beta\gamma}\otimes e\sp
\gamma\otimes b\sp{\alpha\beta} & = &{c\su 1}\sb \alpha{c\su 2} \sb
\beta{}\sco 0
\otimes {c\su 2}\sb \beta{}\sco 1\otimes b\sp{\alpha\beta}\\
& = &  {c\su 1}\sco 0{c\su 2} \sb \beta{}\sco 0
\otimes {c\su 2}\sb \beta{}\sco 1 \otimes b\sp \beta{c\su 1}\sco 1 \\
& = & {c\su 1}\sco 0{c\su 2} \sco 0 \sco 0
\otimes {c\su 2}\sco 0\sco 1 \otimes b{c\su 1}\sco 1{c\su 2}\sco 1 \\
& = & \chi\sb M({c\su 1}\sco 0\otimes\sb M{c\su 2} \sco 0)
\otimes b{c\su 1}\sco 1{c\su 2}\sco 1.
\end{eqnarray*}
{}From the covariance properties of the translation map \cite{Brz:tra}
it then  follows that
$$
{c\su 1}\sb \alpha{c\su 2}\sb{\beta\gamma}\otimes e\sp
\gamma\otimes b\sp{\alpha\beta} =  \chi\sb M(\tau(c\sw 2))
\otimes b(Sc\sw 1)c\sw 3 = 1\otimes c\sw 2\otimes bS(c\sw 1)c\sw 3. $$
This also follows from covariance of $\chi_M$ as intertwining $\Delta^2_R$
projected to $P\tens_M P$ with the tensor product coaction $\Delta^1_R\tens
Ad_R$ on $P\tens H$. Hence Condition 2 may be written as
$$
\overleftarrow{\psi}\sp 2(b\otimes\omega(c)) = \omega(c\sw 2)\otimes
b(Sc\sw 1)c\sw 3
$$
which is equivalent to $\Delta\sb R\sp 2\circ\omega = (\omega\otimes
\id)\circ Ad\sb R$.
\endproof
\begin{ex} For a braided group principal bundle $P(M,B,\Psi)$ in
Example~\ref{braided.ex}, Lemma~\ref{sufficient.lemma} holds. Moreover,
Condition~2 in Proposition~3.5 is equivalent to $Ad_R$-covariance of
$\omega$, where $Ad_R$ is the braided adjoint coaction  as in \cite{Ma:lie}.
\label{braided.connection.ex}
\end{ex}
\proof The braided group adjoint action is studied extensively in
\cite{Ma:lie} as the basis of a theory of braided Lie algebras; we
turn the diagrams up-side-down for the braided adjoint coaction and
its properties (or see earlier works by the second author). Firstly,
$\psi(B\tens M)\subset M\tens B$ is immediate since by properties of
$e=1$, $\Psi(B\otimes M)\subset M\otimes B$. Also clear
is that $\Delta^1_R$ coincides with the given braided coaction of $B$ on $P$
and $\Delta^2_R$ coincides with the braided tensor product coaction on
$P\tens P$. $\Delta\sb R\sp 2$ projects to a coaction on $P\tens_M P$ by
Lemma~\ref{sufficient.lemma}.  We show first that $\chi_M:P\tens_M P\to
P\tens B$ intertwines this coaction with the braided tensor product coaction
$\Delta_R^1\tens Ad_R$. We work with representatives in $P\tens P$ and
use the notation\cite{Ma:introp} as in the proof of
Example~\ref{braided.ex}. Branches $\epsfbox{deltafrag.eps}$ labelled
$\Delta$ are the coproduct of $B$; otherwise they are the given
coaction of $B$ on $P$. Thus,
\[ \epsfbox{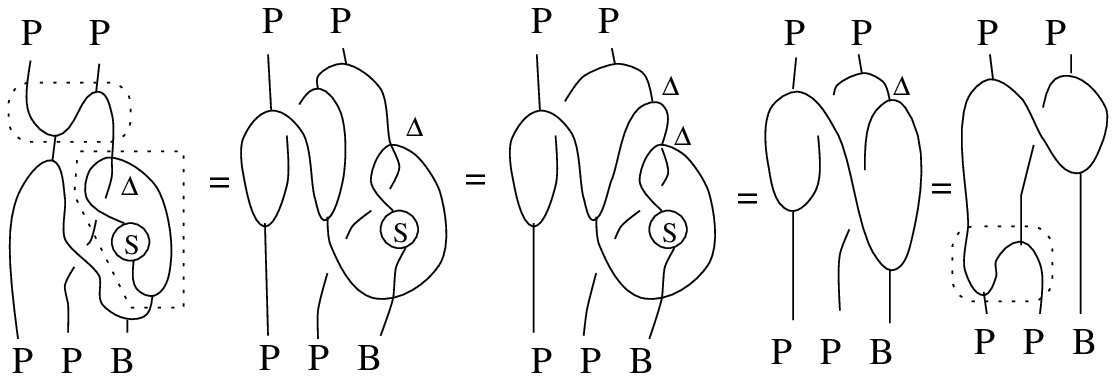}\]
where the upper  box on the left is $\chi_M$ and the lower box is the braided
adjoint coaction $Ad_R$. $S$ denotes the braided antipode of $B$. The tensor
product $\Delta_R^1\tens Ad_R$ uses the braiding and the product of $B$
according to the theory of braided groups\cite{Ma:bg}. The first
equality uses
the homomorphism property of the given coaction of $B$ on $P$. The
second uses
the comodule axiom. The third identifies an `antipode loop' and cancels it
(using associativity and coassociativity, and the braided antipode
axioms). The
fourth equality uses the comodule axiom in reverse and also pushes the
diagram
into the form where we recognise the braided tensor product coaction
$\Delta^2_R$ followed by $\chi_M$. Using this intertwining property of
$\chi_M$, we write the right hand side of Condition~2 in
Proposition~\ref{connection.form.prop} as
\[ \epsfbox{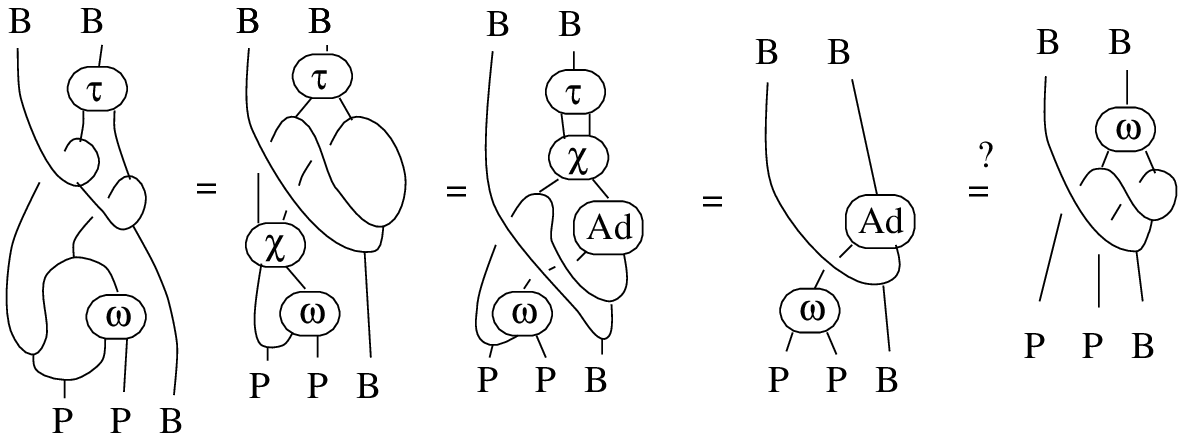}\]
where $\tau=\chi_M^{-1}(1\tens (\ ))$. The left hand side
$\overleftarrow{\psi}\sp 2(b\otimes\omega(c))$ is shown on the right hand side
of the diagram (using associativity of the product in $B$). Hence equality is
equivalent to $\Delta\sb R\sp 2\circ\omega = (\omega\otimes
\id)\circ Ad\sb R$.
\endproof

We remark that in the framework with $C^*$ in place of $C$ as explained in
Proposition~\ref{factor.prop}, we can use for $C^*$ braided groups of
enveloping algebra type, in particular $U({\cal L})$ associated to a
braided-Lie algebra ${\cal L}$ in \cite{Ma:lie} with braided-Lie bracket based
on the properties of the braided adjoint action. In this case one takes
$\omega\in {\cal L}\tens\Omega^1P$ with the corresponding covariance
properties. Using the braided Killing form also in \cite{Ma:lie} one has the
possibility (for the first time) to write down scalar Lagrangians built
functorially from
$\omega$ and its curvature. On the other hand, for a theory of trivial bundles
(in order to have familiar formulae for gauge fields
on the base) one needs to restrict trivialisations and gauge transforms in such
a way that $\omega$ retains its values in $\cal L$. This aspect requires
further work, to be developed elsewhere.

\begin{ex}
Consider $H(M,C, \pi)$,  the $\psi$-principal bundle associated to an
embeddable quantum homogeneous space in
Example~\ref{embeddable.ex}. Assume that $\psi(C\otimes
M)\subset M\otimes C$. Condition~2 in
Proposition~\ref{connection.form.prop} is   equivalent to
\begin{equation}
\Delta\sp 2\sb R\circ\omega\circ\pi = (\omega\otimes
\id)\circ(\pi\otimes\pi)\circ Ad\sb R.
\label{covariance.embeddable}
\end{equation}
In particular, this implies that any inclusion $i: \ker\eps\to H$ such that
$\pi\circ i =\id$ and $\eps(c) = \eps\circ i(c)$
gives rise to the canonical connection 1-form $\omega(c)=(Si(c)\sw 1) \d
i(c)\sw
2$, provided that
$$
(\id\otimes\pi)\circ Ad\sb R \circ i = (i\otimes \id)\circ(\pi\otimes
\pi)\circ Ad\sb R\circ i.
$$
\end{ex}
\proof
In this case $\psi(c\otimes v) = v\sw 1\otimes \pi(u v\sw 2)$,
 and
$\tau(c) = Su\sw 1\otimes\sb M u\sw 2$, for any $c\in C$, $v\in H$ and
$u\in\pi\sp{-1}(c)$. Also
$e = \pi(1)$. The transformation property of $\omega$ now reads
\begin{eqnarray*}
\overleftarrow{\psi}\sp 2(b\otimes\omega(c)) &=&
(Su\sw 1)\sb \alpha u\sw 2 \sb{\beta\gamma}
\omega(\pi(1)\sp \gamma)\otimes \pi(v)\sp{\alpha\beta} \\
&=&  (Su\sw 2)u\sw 3\sb{\beta\gamma}\omega(\pi(1)\sp \gamma\otimes
 \pi(vSu\sw 1)\sp \beta \\
&=&   (Su\sw 2)u\sw 3\sb \gamma\omega(\pi(1)\sp \gamma)\otimes \pi(v(Su\sw 1)
u\sw 4)\\
& = & (Su\sw 2)u\sw 3 \omega(\pi(u\sw 4))\otimes \pi(v(Su\sw 1)u\sw 5) =
\omega(\pi(u\sw 2))\otimes \pi(v(Su\sw 1)u\sw 3) ,
\end{eqnarray*}
where $v\in\pi\sp{-1}(b)$. Choosing $v=1$ we obtain property
(\ref{covariance.embeddable}). The converse is obviously true.
\endproof

Before we describe some concrete examples of connections we construct
connections in  trivial $\psi$-bundles of
Proposition~\ref{trivial.prop}.
\begin{prop} Let $P(M,C,\Phi,\psi,\psi^C,e)$ be a trivial coalgebra
$\psi$-principal bundle such that $\psi(C\otimes M)\subset M\otimes
C$. Let $\beta: C\to \Omega^1M$ be a linear map, $\beta(e)
=0$ and such that
\begin{equation}
\psi^C_{34}\circ\psi_{23}\circ\psi_{12}\circ(\id\otimes\beta\otimes\id)
\circ(\id\otimes\Delta) =
(\beta\otimes\id\otimes\id)\circ(\Delta\otimes\id)\circ\psi^C .
\label{beta.condition}
\end{equation}
Then the map $\omega:\ker\epsilon\to \Omega^1P$,
\begin{equation}
\omega = \Phi^{-1}*\d\Phi + \Phi^{-1}*\beta*\Phi
\label{connection.form.trivial}
\end{equation}
is a connection one-form in $P(M,C,\Phi,\psi,\psi^C,e)$. In particular
for $\beta =0$ we have a trivial connection in
$P(M,C,\Phi,\psi,\psi^C,e)$.
\label{connection.trivial.prop}
\end{prop}
\proof To prove the proposition we will show that $\omega$ satisfies
conditions specified in
Proposition~\ref{connection.form.prop}. Firstly, however, we observe
that the translation map in $P(M,C,\Phi,\psi,\psi^C,e)$ is given by
\begin{equation}
\tau(c) = \Phi^{-1}(c\sw 1)\otimes_M\Phi(c\sw 2).
\label{trans.map.trivial}
\end{equation}
Indeed, a trivial computation shows that $\chi_M(\tau(c)) = 1\otimes
c$, as required. The same computation shows that for any $c\in
\ker\epsilon$,
$$
\chi(\Phi^{-1}(c\sw 1)\d\Phi(c\sw 2) + \Phi^{-1}(c\sw 1)\beta(c\sw
2)\Phi(c\sw 3)) = \chi(\Phi^{-1}(c\sw 1)\otimes\Phi(c\sw 2)) =1\otimes
c,
$$
and therefore Condition 1 of
Proposition~\ref{connection.form.prop} is satisfied by $\omega$.

Now we prove that Condition 2 of
Proposition~\ref{connection.form.prop}  holds for $\Phi^{-1}*\d\Phi$
and $\Phi^{-1}*\beta*\Phi$ separately. For the former the left hand
side of Condition 2 reads
\begin{eqnarray*}
LHS & = &\overleftarrow{\psi}^2(b\otimes \Phi^{-1}(c\sw 1)\otimes\Phi(c\sw 2))
= \Phi^{-1}(c\sw 1)_\alpha\otimes\Phi(c\sw 2)_\beta\otimes
b^{\alpha\beta}\\
& = & \Phi^{-1}(c\sw 1)_\alpha\otimes\Phi(c\sw 2_A)\otimes b^{\alpha
A}
\end{eqnarray*}
On the other hand we use a definition of $\tau$
(\ref{trans.map.trivial})
 and the
properties of $\Phi$ to write the
right hand side of condition 2. as follows:
\begin{eqnarray*}
RHS &= & \Phi^{-1}(c\sw 1)_\alpha\Phi(c\sw 2)_{\beta\gamma}
\Phi^{-1}({e^\gamma}\sw 1)\otimes\Phi({e^\gamma}\sw 2)\otimes
b^{\alpha\beta} \\
& = & \Phi^{-1}(c\sw 1)_\alpha\Phi(c\sw 2)_{\beta\gamma\delta}
\Phi^{-1}({e^\delta})\otimes\Phi({e^\gamma})\otimes
b^{\alpha\beta}\\
& = & \Phi^{-1}(c\sw 1)_\alpha\Phi(c\sw 2_A)_{\gamma\delta}
\Phi^{-1}({e^\delta})\otimes\Phi({e^\gamma})\otimes
b^{\alpha A}\\
& = & \Phi^{-1}(c\sw 1)_\alpha\Phi(c\sw 2_A\sw 1)
\Phi^{-1}(c\sw 2_A\sw 2)\otimes\Phi(c\sw 2_A\sw 3)\otimes
b^{\alpha A}\\
& = &\Phi^{-1}(c\sw 1)_\alpha\otimes\Phi(c\sw 2_A)\otimes b^{\alpha
A} = LHS.
\end{eqnarray*}
To compute the action of $\overleftarrow{\psi}^2$ on the second part
of $\omega$ we will use the shorthand notation
$$
\overleftarrow{\psi}^2(b\otimes\rho) =
\rho_{\underline{\alpha}}\otimes b^{\underline{\alpha}},
$$
for any $b\in C$ and $\rho\in \Omega^1 P$. In this notation Equation
(\ref{beta.condition}) explicitly reads
$$
\beta(c\sw 1)_{\underline{\alpha}} \otimes {c\sw 2}_A \otimes
b^{\underline{\alpha} A} =
\beta({c_A}\sw 1)\otimes {c_A}\sw 2 \otimes b^A,
$$
Using the similar steps as in computation of the action of
$\overleftarrow{\psi}^2$ on the first part
of $\omega$ we find that the right hand side of Condition 2 reads
$$
\Phi^{-1}(c\sw 1)_\alpha\beta(c\sw 2_A\sw 1)\Phi(c\sw 2_A\sw 2)\otimes
b^{\alpha A},
$$
while the left hand side is
$$
\Phi^{-1}(c\sw 1)_\alpha\beta(c\sw 2)_{\underline{\alpha}}\Phi(c\sw 2_A)\otimes
b^{\alpha \underline{\alpha} A} = \Phi^{-1}(c\sw 1)_\alpha\beta(c\sw 2_A\sw
1)\Phi(c\sw 2_A\sw 2)\otimes
b^{\alpha A}.
$$
{}From Proposition~\ref{connection.form.prop} we now deduce that
$\omega$ is a connection one-form as required.
\endproof

Using  similar arguments as in \cite{BrzMa:gau} we can easily show that the
behaviour of
$\beta$  under gauge transformations is exactly the same as in the
case of quantum principal bundles. For example, if
we make a gauge transformation of $\Phi$, $\Phi\mapsto\gamma*\Phi$ and
then view $\omega$ in this new trivialisation then the local
connection one-from $\beta$ will undergo the gauge transformation
\begin{equation}\label{beta.gauge}
\beta\mapsto \gamma^{-1}*\d\gamma +\gamma^{-1}*\beta*\gamma .
\end{equation}

As before, we can specialise this theory to our various sources of
$\psi$-principal bundles. For quantum principal bundles we recover the
formalism in \cite{BrzMa:gau}. For braided principal bundles we make a
computation similar to the one for $\gamma$ in Example~\ref{braided.trivial},
finding that (\ref{beta.condition}) is naturally ensured by requiring that
$\beta:B\to \Omega^1M$ is a morphism in our braided category. Then the same
formulae (\ref{connection.form.trivial}) and transformation law
(\ref{beta.gauge}) etc. apply in the braided case. Indeed, they do not involve
any braiding directly.

Now we construct explicit examples of connections in one of the bundles
described at the end of Section~2.

\begin{ex}
Consider the quantum cylinder bundle
$A\sb q\sp{2|0}[x\sp{-1}](k[x,x\sp{-1}], k[c],\psi, 1)$ in
Example~\ref{cylinder.ex}. Then
$\psi(k[c]\otimes k[x,x\sp{-1}])\subset
k[x,x\sp{-1}]\otimes k[c]$.
The most general connection of the type described in
Proposition~\ref{connection.trivial.prop} has the form
\begin{eqnarray}
\omega (c\sp n)\!\! & = & \sum\sb{k=0}\sp{n-1}(-1)\sp k \qbinom nkq
q\sp{k(k-1)/2} y\sp k\d y\sp{n-k} \nonumber\\
&&\!\!\! + \sum_i \sum_{m=0}^n\sum_{k=0}^m
(-1)^kq^{k((k-1)/2+i)}\qbinom nmq \qbinom mkq \Gamma_{i,m-k}
x^{i}y^k(\d x^{m-k-i})y^{n-m}, \label{connection.braided.line}
\end{eqnarray}
where for all $i\in\bf Z$, $n\in {\bf Z}_{\geq 0}$, $\Gamma_{n,i}\in
k$, $\Gamma_{0,i}=0$.
\label{cylinder.universal.connection.ex}
\end{ex}
\proof  If we set $n=0$ in formula (\ref{cylinder.psi.eqn}) then we
find
$\psi(c\sp l\otimes x\sp m) = q\sp{lm}x\sp m\otimes c\sp l$ and the
first assertion holds. This assertion also follows from
Example~\ref{braided.connection.ex}. We identify $C=k[c]$ by $c_n=c^n$,
as a certain (braided) coalgebra.

It is an easy
exercise to check that a map $\beta : k[c] \to \Omega^1 k[x,x^{-1}]$
satisfies condition (\ref{beta.condition}) if and only if
\begin{equation}
\beta(c^n)
= \sum_i \Gamma_{n,i} x^i\d x^{n-i},
\label{beta.braided.line}
\end{equation}
 where $i\in \bf Z$, $\Gamma_{n,i}\in
k$, $\Gamma_{0,i}=0$.  Now writing the explicit definition of
trivialisation $\Phi$ (\ref{trivialisation.eqn}), and the coproduct of
$c^n$ (\ref{coproduct.braided.line}) we see that $\omega$ in
(\ref{connection.braided.line}) is as in
(\ref{connection.form.trivial}) with $\beta$ given by
(\ref{beta.braided.line}). \endproof

{}From the braided  bundle point of view in Example~\ref{line.ex} on the same
bundle,
we work in the braided category of $\Z$-graded spaces and are allowed for
$\beta$  any degree-preserving that vanishes on $1$. This immediately fixes it
in the form (\ref{beta.braided.line}),  and hence $\omega$ from
(\ref{connection.form.trivial}).

\section{Bundles with general differential structures.}
Let $P(M,C,\psi,e)$ be a $\psi$-principal bundle  as in
Proposition~\ref{bundle.prop}. Let $\CN$ be a
subbimodule of $\Omega\sp 1 P$
such that $\overleftarrow{\psi}\sp 2(C\otimes \CN)\subset \CN\otimes
C$. The map $\overleftarrow{\psi}\sp 2$ induces a map
$\overleftarrow{\psi}\sp 2\sb \CN : C\otimes \Omega\sp 1 P/\CN \to
\Omega\sp 1 P/\CN\otimes C$ and $\CN$  defines a
right-covariant differential structure $\Omega\sp 1(P) = \Omega\sp
1P/\CN$ on $P$. We say that $\Omega\sp 1(P)$ is a differential structure
on $P(M,C,\psi ,e)$.
 \begin{df}
Let $P(M,C,\psi, e)$ be a coalgebra $\psi$-principal bundle and let
$\psi(C\otimes M)\subset M\otimes C$. Assume that $\CN\subset\Omega\sp 1
P$ defines a differential structure $\Omega\sp 1(P)$ on $P(M,C,\psi
,e)$. A connection in $P(M,C,\psi,e)$ is a  left $P$-module
projection $\Pi :\Omega\sp 1 (P)\to \Omega\sp 1(P)$
such that $\ker \Pi = P\Omega\sp 1(M)P$ and $\overleftarrow{\psi}\sp
2\sb \CN(\id\otimes\Pi)  = (\Pi
\otimes \id)\overleftarrow{\psi}\sp 2\sb \CN$.
\label{connection.general.df}
\end{df}

Similarly as for the universal differential calculus case, a
connection in $P(M,C,\psi ,e)$ can be described by its connection
one-form. First we consider the vector space
$\CM = (P\otimes\ker\eps)/\chi(\CN)$ with a canonical surjection $\pi\sb\CM
:P\otimes\ker\eps\to \CM$. Since $\chi$ is a
left $P$-module map, $\chi(\CN)$ is a left $P$-sub-bimodule of
$P\tens\ker\eps$. Therefore $\CM$ is a left
$P$-module and $\pi\sb\CM$ is a left $P$-module map. The action of $P$
on $\CM$ is defined by
$$
u\cdot\upsilon = \sum_i\pi\sb\CM(uv\sb i\otimes c\sp i),
$$
where $u\in P$, $\upsilon\in \CM$ and $\sum\sb iu\sb
i\otimes c\sp i\in \pi\sb\CM\sp{-1}(\upsilon)$. We denote $\CL =
\pi\sb\CM(1\otimes\ker\eps)$. The left $P$-module structure of $\CM$
implies that  for every
element $\upsilon \in\CM$, there exist $u\sb
i\in P$ and $\lambda^i\in\CL$ such that
$\upsilon = \sum\sb i u_i\cdot\lambda^i$. Therefore there is a natural
projection $P\tens\CL\to \CM$.

We assume that $\psi(C\otimes M)\subset M\otimes C$ and
hence the map
$\phi$ can be defined. For any $u\in P$,
$c\in \CM$ and $b\in C$ we have
$$
\phi (b\otimes u\otimes c) = \phi(b\otimes\chi(n)) =
(\chi\otimes \id)\circ\overleftarrow{\psi}\sp 2(b\otimes n) \in
\chi(\CN)\otimes C,
$$
where $n\in \CN$ is such that $\chi(n) = u\otimes c$. We used the
fact that $\overleftarrow{\psi}\sp 2(C\otimes \CN)\subset \CN\otimes C$.

Therefore we can define a map $\phi\sb \CN: C\otimes
\CM \to \CM\otimes
C$ by the diagram
$$
\begin{picture}(290,125)(20,1)
\put(10,110){\makebox(0,0){$C\otimes P\otimes\ker\epsilon$}}
\put(50,110){\vector(1,0){78}}
\put(70,115){$\id\otimes\pi\sb \CM$}
\put(170,110){\makebox(0,0){$C\otimes \CM$}}
\put(205,110){\vector(1,0){60}}
\put(275,110){\makebox(0,0){$0$}}
\put(170,100){\vector(0,-1){50}}
\put(180,85){$\phi\sb \CN$}
\put(170,45){\makebox(0,0){$\CM\otimes C$}}
\put(215,45){\vector(1,0){50}}
\put(275,45){\makebox(0,0){$0$}}
\put(20,100){\vector(0,-1){50}}
\put(30,75){$\phi$}
\put(20,45){\makebox(0,0){$P\otimes \ker\epsilon\otimes C$}}
\put(55,45){\vector(1,0){78}}
\put(87,53){\makebox(0,0){$\pi\sb \CM\otimes \id$}}
\end{picture}
$$
The map $\chi$ induces a map $\chi\sb\CN :\Omega^1(P)\to
\CM$ by the commutative diagram
$$
\begin{picture}(280,155)(20,1)
\put(20,140){\makebox(0,0){$\Omega\sp 1 P$}}
\put(45,140){\vector(1,0){98}}
\put(87,145){$\pi\sb \CN$}
\put(170,140){\makebox(0,0){$\Omega\sp 1 (P)$}}
\put(205,140){\vector(1,0){50}}
\put(265,140){\makebox(0,0){$0$}}
\put(170,130){\vector(0,-1){50}}
\put(180,105){$\chi\sb \CN$}
\put(170,75){\makebox(0,0){$\CM$}}
\put(205,75){\vector(1,0){50}}
\put(265,75){\makebox(0,0){$0$}}
\put(170,65){\vector(0,-1){50}}
\put(170,5){\makebox(0,0){$0$}}
\put(20,130){\vector(0,-1){50}}
\put(30,105){$\chi$}
\put(20,75){\makebox(0,0){$P\otimes \ker\epsilon$}}
\put(50,75){\vector(1,0){90}}
\put(100,83){\makebox(0,0){$\pi\sb \CM$}}
\put(20,65){\vector(0,-1){50}}
\put(20,5){\makebox(0,0){$0$}}
\end{picture}
$$
Clearly, $\chi\sb\CN$ is a left $P$-module map, i.e.,
$\chi_\CN(u\d v) = u\cdot \chi_\CN(\d v)$.
We can use the map $\chi\sb\CN$ to obtain another description of
$\phi\sb \CN$.
\begin{lemma}
The following diagram
$$
\begin{picture}(200,125)(20,1)
 \put(10,100){\makebox(0,0){$C\otimes \Omega\sp 1 (P)$}}
 \put(45,100){\vector(1,0){98}}
\put(87,105){$\overleftarrow{\psi}\sp 2\sb \CN$}
\put(170,100){\makebox(0,0){$\Omega\sp 1 (P)\otimes C$}}
 \put(170,90){\vector(0,-1){50}}
 \put(180,65){$\chi\sb \CN\otimes \id$}
 \put(170,35){\makebox(0,0){$\CM\otimes C$}}
 \put(20,90){\vector(0,-1){50}}
 \put(30,65){$\id\otimes\chi\sb \CN$}
\put(10,35){\makebox(0,0){$C\otimes \CM$}}
 \put(50,35){\vector(1,0){81}}
\put(100,43){\makebox(0,0){$\phi\sb \CN$}}
\end{picture} $$
is commutative.
\label{covariance.chiN.lemma}
\end{lemma}
\proof
We take any $\upsilon\in \Omega\sp 1(P)$, $c\in C$ and $\tilde\upsilon
\in \pi\sb \CN\sp{-1}(\upsilon)$ and compute
\begin{eqnarray*}
\phi\sb \CN\circ(\id\otimes\chi\sb \CN)(c\otimes\upsilon)\!\!\! & = &\!\!\!
\phi\sb
\CN\circ(\id\otimes \pi\sb \CM)(c\otimes \chi(\tilde\upsilon)) =
(\pi\sb \CM\otimes \id)\circ\phi(c\otimes \chi(\tilde\upsilon))\\
& = & (\pi\sb \CM\otimes \id)\circ(\chi\otimes \id)\circ
\overleftarrow{\psi}\sp 2(c\otimes \tilde\upsilon)\\
& = & (\chi\sb \CN\otimes
\id)\circ (\pi\sb \CN\otimes \id)\circ\overleftarrow{\psi}\sp
2(c\otimes\tilde\upsilon)
= (\chi\sb \CN\otimes
\id)\circ\overleftarrow{\psi}\sp
2\sb \CN(c\otimes\upsilon).
\end{eqnarray*}
\endproof

Using arguments similar to the proof of Example~4.11 of \cite{BrzMa:gau}
 and the  definition of a coalgebra
$\psi$-principal bundle
$P(M,C,\psi,e)$ we deduce that
\begin{equation}
0\rightarrow P\Omega\sp 1(M)P\rightarrow \Omega\sp 1(P)
\stackrel{\chi\sb \CN}{\rightarrow}
 \CM\rightarrow 0
\label{exact.sequence.general}
\end{equation}
is a short exact sequence of left $P$-module maps.

\begin{prop}
Connection in $P(M,C,\psi,e)$ with differential structure induced by
$\CN$ is equivalent to a left $P$-module splitting $\sigma\sb \CN$ of the
sequence (\ref{exact.sequence.general}),
such that
$$
(\sigma\sb \CN\otimes \id)\circ\phi\sb{\CN} =
\overleftarrow{\psi}\sp 2\sb{\CN}\circ(\id\otimes\sigma\sb \CN) .
$$
\end{prop}
\proof
We use Lemma~\ref{covariance.chiN.lemma} to deduce the covariance
properties of $\chi\sb \CN$ and then preform calculation similar to the
proof of Proposition~\ref{proposition.splitting.universal}.
\endproof

To each connection $\Pi$ we can associate its connection one form
$\omega:\CL\to\Omega\sp 1(P)$ by $\omega(\lambda) =
\sigma\sb\CN(\lambda)$.

Similarly to the case of universal differential structure, one proves
\begin{prop}
Let $\Pi$ be a connection in $P(M,C,\psi,e)$ with differential
structure defined by $\CN\subset \Omega\sp 1 P$. Then, for all
$\lambda\in\CL$ the  connection
1-form $\omega :\CL\to\Omega\sp 1 (P)$ has the following
properties:

1. $\chi_\CN\circ\omega(\lambda) = \lambda$,

2. For any $b\in C$, $\overleftarrow{\psi}\sp 2_\CN (b\otimes\omega(\lambda)) =
 \tilde c\su 1\sb \alpha \tilde c\su 2\sb{\beta\delta}  \omega(\pi\sb
\CM(1\otimes e\sp
\delta))\otimes b\sp{\alpha\beta}, $
where $\tilde c\su 1\otimes\sb M\tilde c\su 2$ denotes the translation
map $\chi\sb M\sp{-1}(1\otimes \tilde c)$, and $\tilde c \in\ker\eps$
is such that $\pi_\CM(1\otimes\tilde{c}) = \lambda$.

Conversely, if $\CM$ is isomorphic to $P\tens\CL$ as a left $P$-module
and $\omega$ is any linear map $\omega :\CL\to
\Omega\sp 1(P)$ obeying conditions 1-2, then there is a unique
connection $\Pi = \cdot\circ(\id\tens\omega)\circ\chi_\CN$ in
$P(M,C,\psi, e)$ such
that $\omega$ is its connection 1-form.
\label{connection.form.general.prop}
\end{prop}

In the setting of \cite{BrzMa:gau} the condition $P\tens\CL = \CM$ is
always satisfied for quantum principal bundles, and $\CL =
\ker\eps/{\cal Q}$, where ${\cal Q}$ is an $Ad_R$-invariant right ideal in
$\ker\eps$ that generates the bicovariant differential structure on
the structure quantum group $H$ as in \cite{Wor:dif}. The detailed
analysis of braided group principal
bundles with general differential structures will be presented
elsewhere. Here we remark only that it seems natural to assume that
$\CM = P\und\tens \CL$ and then choose $\CL$
to be the space dual to the braided Lie algebra $\cal L$ as discussed
in Section~3. This choice of $\CL$ is justified by the fact that from
the properties of the maps $\phi$ and $\phi_\CN$ it follows
  that the space $\CL$ is invariant under  the braided adjoint
coaction (cf. Example~\ref{braided.connection.ex}).

We complete this section with an explicit example of differential
structures and connections on the quantum cylinder bundle in
Example~\ref{cylinder.ex}
(cf. Example~\ref{cylinder.universal.connection.ex}).
\begin{ex}
\rm We consider the quantum cylinder bundle of
Example~\ref{cylinder.ex} (cf. Example~\ref{line.ex}) and we work with
differential structures on $A_q^{2|0}$ classified  in
\cite{BrzDab:dif}. Using definition of $\psi$ (\ref{cylinder.psi.eqn})
one easily finds that there are two differential structures for which
the covariance condition $\overrightarrow{\psi}^2 (k[c]\otimes \CN
)\subset \CN\tens k[c]$ is satisfied. The subbimodules $\CN$ are generated
by
$$
(1+s)x\tens x- x^2\tens 1 - 1\tens x^2, \quad y\tens x - qxy\tens 1 -
q1\tens xy +qx\tens y, \quad (1+q)y\otimes y -y^2\otimes 1 -1\otimes
y^2,
$$
where $s\in  k$ is a free parameter, in the first case, and by
$$
(1+q)x\tens x- x^2\tens 1 - 1\tens x^2, \quad y\tens x - xy\tens 1 -
q1\tens xy +x\tens y, \quad (1+q)y\otimes y -y^2\otimes 1 -1\otimes
y^2,
$$
in the second case. In both cases the modules of 1-forms
$\Omega^1(A^{2|0}_q)$ are generated by the exact one-forms $\d x$ and
$\d y$. Definitions of the $\CN$ imply the following relations in
$\Omega^1(A^{2|0}_q)$
$$
x\d x = s\d x x,\quad x\d y = q^{-1}\d yx, \quad y\d x =q\d xy, \quad
y\d y = q\d yy,
$$
in the first case, and
$$
x\d x = q\d x x,\quad x\d y = \d yx, \quad y\d x =q\d xy +(q-1)\d xy, \quad
y\d y = q\d yy,
$$
in the second one. In both cases $(A^{2|0}_q[x^{-1}]
\otimes\ker\eps)/\chi(\CN) = A^{2|0}_q[x^{-1}]\otimes \CL$, where $\CL$ is a
one-dimensional vector space spanned by $\lambda = \pi\sb\CM(1\otimes
c)$ and can be therefore identified with a subspace of $k[c]$ spanned
by $c$. Also in both cases the most general connection is given
by
$$
\Pi(\d x) = 0, \qquad \Pi(\d y) =\d y +\alpha \d x,
$$
where $\alpha\in k$, and extended to the whole of
$\Omega^1(A^{2|0}_q[x^{-1}])$ as a left $A^{2|0}_q[x^{-1}]$-module
map. The corresponding connection  one form reads
$$
\omega(\lambda) = \d y +\alpha \d x .
$$
The bundle is trivial and this connection can be described by the
map $\beta :k[c] \to \Omega(k[x,x^{-1}])$ as in
Proposition~\ref{connection.trivial.prop}
(cf. Eq.~(\ref{beta.braided.line})) with $\beta(c^n) =0$ if $n\neq 1$
and
$\beta(c) = \alpha\d x$.
\end{ex}

\end{document}